   \def\obs#1{{\bf (*** #1 ***)} }

 \def\obs#1{}     

\documentclass[twoside,letterpaper,draft,11pt]{amsart}

\usepackage{amsmath}               
\usepackage{amsthm}                
\usepackage{latexsym}

\usepackage{xspace}
\usepackage{amscd}

\title[Globalization of  twisted partial actions]{Globalization of twisted partial actions}

\author[M.\ Dokuchaev]{M. Dokuchaev}
\address{Instituto de
Matem\'atica e Estat\'\i stica,
Universidade de S\~ao Paulo, 
05508-090 S\~ao Paulo, SP, Brasil}
\email{dokucha@ime.usp.br}
\author[R.\ Exel]{R. Exel}
\address{Departamento de Matem\'{a}tica,
 Universidade Federal de Santa Catarina,
 88040-900 Florian\'opolis, SC,  Brasil}
\email{exel@mtm.ufsc.br}
\author[J.\ J.\ Sim{\'o}n]{J. J. Sim{\'o}n}
\address{Departamento de Matem\'{a}ticas,
 Universidad de Murcia,
 30071 Murcia, Espa\~{n}a}
\email{jsimon@um.es}
\thanks{{This work was partially supported by CNPq of Brazil and Secretar\'{\i}a de Estado de Universidades e Investigaci\' on del MEC, Espa\~ na.\\{\bf 2000 Mathematics Subject Classification}:
Primary 16W50; Secondary   16S35, 16W22.\\
{\bf Key words and phrases:} partial action, twisting, crossed
product, globalization.}}

\newtheorem{teo}{Theorem}[section]
\newtheorem{defi}[teo]{Definition}
\newtheorem{lema}[teo]{Lemma}

\newtheorem{prop}[teo]{Proposition}

\theoremstyle{remark}
\theoremstyle{remark}

\newcommand{\0}{\theta}

\newcommand{\D}{{\mathcal D}}

\newcommand{\de}{\delta}

\newcommand{\af}{\alpha}

\newcommand{\bt}{\beta}

\newcommand{\lb}{\lambda}

\newcommand{\Lb}{\Lambda}

\newcommand{\p}{{\bf Proof. }}

\newcommand{\fim}{\hfill\mbox{$\Box$}}

\newcommand{\A}{{\mathcal A}}
\newcommand{\R}{{\mathcal R}}
\newcommand{\B}{{\mathcal B}}

\newcommand{\M}{{\mathcal M}}

\newcommand{\U}{{\mathcal U}}

\newcommand{\F}{{\mathcal F}}

\newcommand{\m}{{}^{-1}}
\newcommand{\mt}{\mapsto}
\def\f{\varphi}
\def\e{\varepsilon}

\def\ndv{\ {\mid \kern -0.7 em {\scriptstyle \not}} \ \ }

\def\nd{\ {\mid \kern -0.4 em {\scriptstyle \not}} \ \ }

\begin{document}

\date{\today}

\begin{abstract} Let $\A $ be a unital ring which is a product of
  possibly infinitely many 
  indecomposable rings. We establish
a criteria for the existence of a globalization for a given
twisted partial action of a group on $\A .$ If the globalization
exists, it is unique up to a certain equivalence relation  and,
moreover, the crossed product corresponding to the twisted partial
action is Morita equivalent to that corresponding to its
globalization. For arbitrary unital  rings the globalization
problem is reduced to an extendibility property of the multipliers
involved in the twisted partial action.
\end{abstract}

\maketitle

\begin{section}{Introduction}\label{sec:intro}

The relationship between partial isomorphisms and global ones
proved to be relevant in diverse areas of mathematics, such as
operator theory, topology, logic, graph theory,  differential
geometry, group theory (with particular importance for geometrical
and combinatorial group theories) and semigroup theory  (see
\cite{KL}, \cite{L}). Moreover, it na\-tu\-ral\-ly can be found as
a significant feature even in the basics of Galois Theory of
fields, or in that of inner product spaces in the form of Witt
Lemma (see \cite{L}). The latter was the origin to the term ``Witt
property'', naming an important situation in which every partial
isomorphism lies beneath a global one. Other remarkable related
properties which reflect certain ``globalization'' phenomenons are
the so-called ``HNN-property'' and that of ``homogeneity'' (see \cite{L}).\\

When a group $G$ is given and a partial isomorphism of a structure
is attached to each element of $G$ such  that the composition of
partial isomorphisms respects the group operation, we speak about
a partial action of $G.$ The  formal definition of this concept
was given in \cite{E1}, and the study of globalizations (also
called enveloping actions) of partial group actions was initiated
in the PhD Thesis of F. Abadie \cite{Abadie} (see also
\cite{AbadieTwo}) and independently by J. Kellendonk and M. Lawson
in \cite{KL}. Further results on globalizations of partial actions
of groups were obtained in \cite{S}, \cite{DE}, \cite{F} and
\cite{DdRS}. In particular, in \cite{DE}  a criteria for the
existence of a globalization for a partial group action on a
unital ring was given, which was further generalized for partial
group actions on left $s$-unital rings in \cite{DdRS},
incorporating a new ingredient to the criteria, which is essential
for the $s$-unital context. One should notice that  the concept of
a globalization of a partial action of a group on a ring $\A ,$ as
treated in  \cite{AbadieTwo}, \cite{DE} and \cite{DdRS} (see also
Definition~\ref{twistdef} below), has a rigid condition that  $\A
$ is an ideal in the ring under the global action. This is imposed
in the definition because of the way 
  the restriction works:
any global group action on a ring canonically restricts to a
partial action on a two-sided ideal. However, it may happen that a
given partial action restricts to a subring which is not an ideal,
and dropping the considered restriction from the definition one
comes to the so-called weak globalizations, on
which interesting results have been obtained in \cite{F}.\\

Twisted partial actions of locally compact groups on
$C^*$-algebras were introduced in \cite{E0} serving the general
construction of $C^*$-crossed products. The power of the notion
allows
  one 
  to prove (see \cite[Theorem 7.3]{E0}) that any second
countable $C^*$-algebraic bundle with stable unit fiber can be
obtained as a result of the construction. Its generality permits
to characterize important classes of $C^*$-algebras as
$C^*$-crossed products (see \cite{ELQ}). Twisted partial actions
of groups on abstract rings and corresponding crossed products
were introduced in \cite{DES1}, in which an algebraic analog of
the above mentioned stabilization theorem was also proved. As a
recent application of the construction, it was  proved in
\cite{E2}, that given a field $F$ with ${\rm char} F =0$ and a
subnormal subgroup $H \lhd N \lhd G,$ there is a twisted partial
action of the factor group $G/N$ on the group algebra $F [N/H]$
such that the Hecke  algebra ${\mathcal H}(G,H)$ is
isomorphic to the partial crossed product $F [N/H] \ast G/N.$
 Separability, semisimplicity and Frobenius properties of crossed products by twisted partial actions were investigated in 
 \cite{BLP}.\\

The globalization problem for partial actions on rings is
essential for the further study of general crossed products and
related topics, in particular,  the above mentioned result from
\cite{DE} 
  turned out to be useful
  in the Galois Theory of
partial actions \cite{DFP}, in M. Ferrero's result \cite{F}, and
in the series of recent preprints by D. Bagio, W. Cortes, M.
Ferrero, J. Lazzarin, H. Marubayashi and A. Paques  \cite{BCFP}, \cite{C},
\cite{CF}, \cite{CFM}, \cite{FL}.  Thus a solution for the globalization problem
for the twisted context seems to be rather desirable, which is our
present purpose.\\

Our paper is organized as follows. After giving some preliminaries
in Section~\ref{sec:preliminaries},  we point out in
Section~\ref{sec:morita} a Morita equivalence result with respect
to the partial and global crossed products, which is an immediate
consequence  of the arguments given in \cite{DE} and \cite{DdRS}.
In Section~\ref{sec:reducing}, we reduce the question of the
existence of a globalization of a twisted partial action of a
group on a  unital ring to an extendability property of the
  multipliers involved. 
Next we pass to consider twisted partial
actions on rings $\A $ which are (finite or infinite) products of
indecomposable unital rings, first giving in
Section~\ref{sec:transitive} some preliminary facts on transitive
twisted partial actions. In Section~\ref{sec:goingto} we construct
a more manageable twisted partial action, equivalent to the given
one. The idea is borrowed from the corestriction map
 in Homological Algebra, which turns out to be  adaptable to the
 non-commutative partial framework. In Section~\ref{sec:existence}
we establish our main existence result whereas  in
Section~\ref{sec:uniqueness} we prove that two globalizations of a
given twisted partial action of a group on $\A $ are equivalent in
  some natural sense from the homological point of view, 
provided
that the rings under the global actions are unital.

\end{section}

\begin{section}{Some preliminaries}\label{sec:preliminaries}

 We remind the reader that the multiplier ring of $\M(\A)$ of an associative
non-necessarily unital ring $\A$ is the
 set $$\M(\A)= \{(R,L) \in {\rm End}(_{\A} \A) \times {\rm End}(\A_{\A}) :
(aR)b = a(Lb) \, \mbox {for all}\, a,b
 \in \A \}$$ with component-wise addition and multiplication (see
\cite{DE}  or \cite[3.12.2]{Fillmore} for more details). Here we
use the right hand side notation   for homomorphisms of left
$\A$-modules, whereas for homomorphisms of right modules the usual
notation is used. Thus given $R : _{\A}\A  \to _{\A}\A,$  $L :
{\A}_{\A} \to {\A}_{\A}$ and $ a \in \A$ we write $a \mapsto
 a R $ and   $a \mapsto L a. $ For a multiplier $w = (R,L) \in \M(\A)$ and an element
 $a \in \A$ we set $a w = a
R$ and $w a = L a,$ so that  one always has $(a w ) b = a (w b)$
$(a,b \in \A).$ Given a ring isomorphism $\phi : \A \to {\A }',$
the map $\M(\A) \ni w \mt \phi \, w \, {\phi }\m \in \M({\A }')$
where $ \phi \, w \, {\phi }\m = ({\phi }\m  \, R \, {\phi } ,
\phi \, L \, {\phi }\m ), w =(R, L),$
is an isomorphism of rings.\\

\begin{defi}\label{twistdef} A twisted partial action of a group  $G$
 on  $\A$ is a triple
$$\af = (\{\D_g \}_ {g \in G}, \{{\af}_g \}_ {g \in G},
 \{w_{g,h} \}_ {(g,h) \in G\times G}),$$ where for each $g \in G,$ ${\D}_g$ is
a two-sided ideal in $\A,$ ${\af }_g$ is an isomorphism of rings
${\D}_{g\m} \to  {\D}_g,$ and for each $(g,h) \in G\times G,$ $
w_{g,h}$ is an invertible element from $\M(\D_g \cdot \D_{gh}),$
satisfying the following
postulates, for all $g,h$ and $t$ in $G:$\\

 (i)  $\D_g^2 = \D_g, \; \D_g \cdot \D_h = \D_h \cdot \D_g;$\\

 (ii) $\D_1 = \A$ and ${\af}_1$ is the identity map of $\A;$\\

 (iii) ${\af}_g({\D}_{g\m} \cdot {\D}_{h}) = {\D}_{g} \cdot {\D}_{gh};$\\

 (iv) ${\af}_g \circ {\af}_h (a) = w_{g,h} {\af}_{gh}(a) w\m_{g,h}, \; \forall
a \in  {\D}_{h\m} \cdot {\D}_{h\m g\m};$\\

 (v) $w_{1,g} = w_{g,1} = 1;$\\

 (vi) $\af_g(a w_{h,t}) w_{g, ht} = {\af}_g(a) w_{g,h} w_{gh,t}, \;
\forall a \in \D_{g\m} \cdot \D_h \cdot \D_{ht}.$

\end{defi}

 As it is commented in \cite{DES1}, it
follows from  (i) that  a finite  product of ideals $\D_g \cdot
\D_h \cdot \ldots$ is  idempotent, and $${\af } _g (\D_{g\m} \cdot
\D_h \cdot \D_f) = \D_{g} \cdot \D_{gh} \cdot \D_{gf}$$ for all
$g,h,f \in G,$ by (iii). Thus all multipliers in (vi) are
applicable.\\

We say that $\af $ is {\it global} if $\D_g = \A$ for all $g \in G.$
Observe that given a twisted global action
$${\bt } = ({\B }, \{{\bt}_{x} \}_ {x \in G},
 \{u_{x,y} \}_ {(x,y) \in G\times G})$$ of a group $G$ on a (non-necessarily unital) ring
 $\B ,$ one may restrict $\bt $ to a two-sided ideal $\A $ of $\B
 $ such that $\A $ has $1_{\A }$ as follows. Setting ${\D }_x
 = \A \cap \bt _x (\A ) =  \A \cdot \bt _x (\A ) $ we see that each
 $\D _x$ has $1$ which is $1_{\A } \bt _x (1_{\A }).$ Then putting
 $\af _x = \bt _x \mid _{\D _{x\m }},$ the items
(i), (ii), (iii) of Definition~\ref{twistdef}
 are satisfied. Furthermore, defining $w_{x,y} = u_{x,y}1_{\A } \bt _x (1_{\A }) \bt _{xy} (1_{\A
 })$ we have that (iv), (v) and (vi) are also satisfied, so  we
 indeed have obtained a twisted  partial action on $\A .$

\begin{defi}\label{globalization} A twisted global  action
$$\bt = ({\B },  \{{\bt}_g \}_ {g \in G},
 \{u_{g,h} \}_ {(g,h) \in G\times G})$$
of a group $G$ on an associative (non-necessarily unital) ring
$\B$ is said to be a glo\-ba\-li\-za\-tion (or an enveloping
action) for the partial action $\af$ of $G$ on $\A$
if there exists a monomorphism $\f: \A \to \B$ such that:\\

\noindent  {\it (i)} $\f(\A)$ is an ideal in $\B,$\\

\noindent {\it (ii)} $\B = \sum_{g \in G} \bt_g(\f(\A)),$\\

\noindent  {\it (iii)}  $\f(\D_g) = \f(\A) \cap {\bt}_g(\f(\A))$ for any $g\in G,$\\

\noindent {\it (iv)} $ \f \circ \af_g   = \bt_g \circ \f $ on $\D_{g\m}$ for any $g\in G,$\\

\noindent {\it (v)} $ \f (a w_{g,h}) = \f (a) u_{g,h},$ $  \f (
w_{g,h} a) =  u_{g,h} \f (a)$ for any $g, h \in G$ and $a \in \D_g
 \D_{gh}.$\\

\noindent If this is the case,  we shall say that $\af$ is
globalizable.
\end{defi}

Our {\it rings} will be always associative,  non-necessarily with
unity in general.  Given a unital ring $\A ,$ we shall denote by
$\U ( \A )$ the group of invertible elements of $\A ,$ whereas
$1_{\A }$ will stand for the unity of $\A .$ If $\bt $ is a
globalization for $\af ,$ then $\B $ may have no unity even when
$\A $ is a ring with $1_{\A},$ if, however, $\B $ has $1_{\B },$
then we say that $\bt $ is a {\it unital } globalization and $\af
$ is {\it unitally globalizable}. The multipliers $w_{x,y}$ of a twisted partial  action will be
also denoted by $w[x,y].$\\

\end{section}

\begin{section}{Morita equivalence}\label{sec:morita}

It was proved in \cite{DE} that if $\af $ is  a (non-twisted)
partial action of a group $G$ on a unital ring $\A $ and $\bt $ is
a globalization of $\af $ with $G$ acting (globally) on $\B $ such
that $\B $ has $1_{\B } $ then the partial skew group ring $\A
\ast _{\af } G $ is Morita equivalent to the (global) skew group
ring $\B \ast _ {\bt} G.$ The proof is easily adaptable for the
left $s$-unital case, as it was shown   in \cite{DdRS}. We recall
that a ring $\A $ is called left $s$-unital if for any $a \in \A$
there exists an element $e\in \A$ such that $ea = a.$
Equivalently, for any finitely many $a_1, \ldots , a_k  \in \A $
there exists an element $e \in \A $ with $ea_i = a_i, i=1, \ldots
k$ (see \cite[Lemma 2.4]{DdRS}). Any left $s$-unital ring is
obviously idempotent, and for them the  Morita theory of
idempotent rings developed in  \cite{GS} is applicable. We also
recall that given a twisted partial action $\af$  of  $G$ on $\A,$
the crossed product $\A \ast _{\af} G$ is the direct sum:
$$\bigoplus_{g \in G} \D_g \de_g,$$ in which the $\de_g$'s are symbols,  and the
multiplication is defined by the rule: $$(a_g {\de}_g) \cdot (b_h
{\de}_h) = {\af}_g ( {\af}\m_{g}(a_g)b_h) w_{g,h} {\de}_{gh}.$$
\noindent Here $w_{g,h}$ acts  as a right multiplier on ${\af }_g
( {\af }\m_{g}(a_g)b_h) \in {\af }_g(\D_{g\m}\cdot \D_h) = \D_g
\cdot \D_{gh}.$ The associativity of this construction was
established in  \cite{DES1}.\\

We point out the next:

\begin{teo}\label{morita} Let $\af $ be
a globalizable twisted partial action of $G$ on a left $s$-unital
ring $\A $ and let $\bt $ be a globalization of $\af $ with $G$
acting globally on $\B .$ Then the crossed products $\A \ast _
{\af} G$ and  $\B \ast _ {\bt} G$ are Morita equivalent.
\end{teo}
 \p We  observe that the Morita context and the
 arguments given in \cite{DE} with the adaptation for the
 $s$-unital case in \cite{DdRS} work nicely, and we give only a sketch of the proof. 
 
  Since  $\A $ is left
 $s$-unital,  so too is $ \B ,$ as a sum of left $s$-unital ideals (see \cite[Remark 2.5]{DdRS}). By \cite[Theorem 3.1]{DdRS} each $\D _g$ is left $s$-unital and  consequently, $\R= \A \ast _
{\af} G$ and  $\R'= \B \ast _ {\bt} G$ are also left $s$-unital. We  view  $\R$ as a subring of $\R'.$ Consider   $M,N\subseteq \B\ast_\bt G$ given by
$$ M = \left\{\sum_{g\in G}c_g u_g : c_g\in \A \right\} \;\;\; \mbox{and} \;\;\;  N = \left\{\sum_{g\in G}c_g
u_g : c_g\in \bt_g(\A) \right\}.$$ The proofs of Propositions
5.1 and 5.2 of \cite{DE} do not use the fact that the  ring has $1$ and work without essential changes in the twisted context. According to these propositions $M$ is an $\R$-$\R'$ - bimodule and $N$ is a $\R'$-$\R$ - bimodule. While checking this, one should keep  in mind that 
all ideals $\bt _x (\A )$ $(x \in G)$ are idempotent,
and therefore, $u[y,z] \bt _x (\A ) = \bt _x (\A ) u[y,z] = \bt _x
(\A )$ for any $x,y,z \in G,$ as invertible multipliers preserve
idempotent ideals.

Next, let $\tau: M \otimes_{\R'} N \to \R'$ and  $\tau': N \otimes_{\R} M \to \R'$ be given by $\tau(m\otimes n) = mn$ and  $\tau'(n\otimes m) = nm.$ Obviously $\tau$ is an $\R$-bimodule map and $\tau'$ is an $\R'$-bimodule map. Similarly as in  the proof of Theorem 4.1 of \cite{DdRS} one directly verifies that $\tau(M\otimes_{\R'}N) = \R$ and   $\tau'(N\otimes_{\R}M) = \R'.$ 

The notion
of the Morita context for idempotent rings  is the same as for rings with $1$ (see \cite{GS}), thus we have a Morita context $(\R, \R',M, N,\tau,\tau')$ with idempotent $R$ and $R',$  and surjective $\tau$ and $\tau',$ therefore by \cite{GS}, the rings $\R$ and ${\R}'$ are Morita equivalent. \fim \\

\end{section}

\begin{section}{Reducing the globalization problem to an extendability property
of multipliers}\label{sec:reducing}

If $\af $ is a globalizable twisted partial action of $G$ on a
unital ring  $\A ,$ then  each $\D _g$ $(g\in G)$  is a unital
ring, as $\f(\D_g) = \f(\A) \cap {\bt}_g(\f(\A)).$  In the
non-twisted case the converse is true by \cite[Theorem 4.5]{DE}.
In the  twisted case, without imposing restrictions on $\A ,$  we
are able to show in Theorem~\ref{globreduced} below that the
converse is true provided that  an extendability property of the
partial twisting ${w}$ is verified (the latter is also a necessary
condition).  if each $\D _g $ is a unital ring whose
unity element is denoted by $1_g$, then $1_g$ is a central
idempotent of $\A$ such that $\D_g = 1_g \A,$ and for every $g, h
\in G$  the ring $\D _g \cap \D _h = \D_g \cdot \D_{h}$ is unital
with unity $1_g 1_{h}.$ Consequently $\M(\D_g \D_{gh}) \cong \D_g
 \D_{gh},$ so that ${w}_{g,h}$ is an invertible element of $\D_g
 \D_{gh}$ for any $g, h \in G.$ Moreover, by (iii) of
 Definition~\ref{twistdef},

\begin{equation}\label{triviality}
\af_g (1_{g\m} 1_h) = 1_g 1_{gh},
 \end{equation}
\noindent for any $g, h\in G.$



\begin{teo}\label{globreduced} Let $\af $ be a twisted partial action of a group $G$ on a unital ring
$\A$ such that each $\D_g$ $(g \in G)$ is a unital ring with unity
element $1_g.$ Then $\af $ admits a globalization if and only if
for each pair $(g,h) \in G \times G$ there exists an invertible
element $\widetilde {w}_{g,h} \in {\U}(\A),$  such that
  $\widetilde {w}_{g,h} 1_g 1_{gh} = {w}_{g,h}, $  and
\begin{equation}\label{extend} \af_g( \widetilde{w}_{h,t}\, 1_{g\m})\, \widetilde{w}_{g, ht}
=1_g \, \widetilde{w}_{g,h} \, \widetilde{w}_{gh,t},
\end{equation} for any $g, h, t \in G.$
\end{teo}

\p The ``only if''  part is obvious by taking
$\widetilde{w}_{g,h}={u}_{g,h} \cdot 1_{\A}$ for any $g,h \in G.$\\

For the ``if'' part let $\F = \F(G,\A)$ be the Cartesian product
of the copies of $\A$ indexed by the elements of $G,$ that is, the
ring of all functions  from $G$ into $\A.$ For convenience of
notation $f(g)$ will be also written as  $f|_g$  $(f \in \F, g \in
G).$

For $g \in G$ and $f \in \F$ define $\bt_g(f) \in \F$ by the
formula:
$$\bt_g (f)|_h = \widetilde {w}_{h\m,g} f(g\m h) \widetilde {w}\m_{h\m,g}, \; \; h \in G.$$
Obviously $f \mapsto \bt_g(f)$ determines an automorphism $\bt_g$
of $\F .$ Define $u_{g,h} \in {\U} (\F)$ by setting
$$u_{g,h}|_t =  \widetilde {w}_{t\m,g} \widetilde {w}_{t\m g, h} \widetilde
{w}\m_{t\m,gh}, \;\; g,h,t \in G.$$

For arbitrary $g, h, t, x \in G$ one has that

\begin{align*} & (\bt _g (u_{h,t}) u_{g,ht}) |_x = \widetilde
{w}_{x\m,g} \, u_{h,t} |_{g\m x} \, {\widetilde {w}}\m_{x\m,g}
u_{g,ht} |_x = \\ & ({\widetilde {w}}_{x\m,g} \,  {\widetilde {w}}_{x\m
g,h} \, {\widetilde {w}}_{x\m gh,t} \,  {\widetilde {w}}\m_{x\m g,ht} \,
{\widetilde {w}}\m_{x\m,g}) \,  ({\widetilde {w}}_{x\m,g} \, {\widetilde
{w}}_{x\m g,ht} \, {\widetilde {w}}\m_{x\m,ght}) =\\ & {\widetilde
{w}}_{x\m,g} \, {\widetilde {w}}_{x\m g,h} \, {\widetilde {w}}_{x\m gh,t}
\,  {\widetilde {w}}\m_{x\m,ght},
\end{align*} and
\begin{align*} & (u_{g,h} u_{gh,t}) |_x = ({\widetilde {w}}_{x\m,g} \,
{\widetilde {w}}_{x\m g,h} \, {\widetilde {w}}\m_{x\m , gh} ) \, ({\widetilde
{w}}_{x\m , gh}  \,  {\widetilde {w}}_{x\m gh,t} \,  {\widetilde
{w}}\m_{x\m,ght}) = \\ &{\widetilde {w}}_{x\m,g} \, {\widetilde {w}}_{x\m
g,h} \, {\widetilde {w}}_{x\m gh,t} \,  {\widetilde {w}}\m_{x\m,ght},
\end{align*} which shows that $$\bt _g (u_{h,t}) u_{g,ht} =  u_{g,h} u_{gh,t},
$$i.e. $u$ satisfies the 2-cocycle equality (vi) of
Definition~\ref{twistdef}.

Next for any $f \in \F$ and $g,h, x \in G$ we compute

\begin{align*}  &(\bt_g \circ \bt _h (f) ) |_x = {\widetilde
{w}}_{x\m,g} \, \bt _h (f) |_{g\m x}  \, {\widetilde {w}}\m_{x\m,g} =
{\widetilde {w}}_{x\m,g} \, {\widetilde {w}}_{x\m g,h} \, f(h\m g\m x) \,
{\widetilde {w}}\m_{x\m g,h} \, {\widetilde {w}}\m_{x\m,g}.
\end{align*} Furthermore,

\begin{align*} & (u_{g,h}\, \bt_{gh}  (f) \,  u\m_{g,h} ) |_x = \\ &
({\widetilde {w}}_{x\m,g} \, {\widetilde {w}}_{x\m g,h} \,  {\widetilde
{w}}\m_{x\m,gh}) \, {\widetilde {w}}_{x\m,gh} \, f(h\m g\m x) \,
{\widetilde {w}}\m_{x\m,gh} \, ( {\widetilde {w}}_{x\m,g} \, {\widetilde
{w}}_{x\m g,h} \, {\widetilde {w}}\m_{x\m,gh})\m=
\\ & {\widetilde {w}}_{x\m,g} \, {\widetilde {w}}_{x\m g,h} \, f(h\m g\m
x) \,   {\widetilde {w}}\m_{x\m g,h} \, {\widetilde {w}}\m_{x\m,g},
\end{align*} and consequently,

\begin{equation}\label{betacomp}
\bt_g \circ \bt _h (f) = u_{g,h}\, \bt_{gh}  (f) \, u\m_{g,h},
\end{equation} i.e. $u$ satisfies (iv) of Definition~\ref{twistdef}. Since
trivially $u_{1,g} = u_{g,1} = 1$ for all $g \in G,$ we have a
twisted global action of $G$ on $\F .$

For any $a \in \A$ the element $a 1_g$  belongs to $\D_g$ and the
formula $$\f(a)|_g = \af_{g\m}(a 1_g), \; g \in G,$$ defines a
monomorphism $\f : \A \to \F.$

 Let $\B = \sum_{g \in G} \bt_g(\f(\A)),\; (g \in G).$
Our purpose  is to show  that the restriction of $\bt$ to $\B$ is
a globalization for $\af$. Denote this restriction by the same
symbol $\bt.$ We proceed by checking property (iv) of
Definition~\ref{globalization}.

For $g,h \in G$ and $a \in \D_{g^{-1}}$ we have 
\begin{align*} & \bt_g(\f(a))|_h
=
 \widetilde {w}_{h\m,g}  \f (a)|_{g\m h}  {\widetilde {w}}\m_{h\m,g} ={\widetilde {w}}_{h\m,g}
 \af_{h\m g} (a 1_{g\m h}) {\widetilde {w}}\m_{h\m,g}=\\
&  {w}_{h\m,g}
 \af_{h\m g} (a 1_{g\m h}) 1_{h\m} 1_{h\m g}  {w}\m_{h\m,g}= {w}_{h\m,g}
 \af_{h\m g} (a 1_{g\m} 1_{g\m h}) {w}\m_{h\m,g},
 \end{align*}   in view of  (\ref{triviality}) and the equalities $a= a 1_{g\m},$ ${\widetilde {w}}_{h\m , g} 1_{h\m}1_{h\m g}= w_{h\m , g}.$
 On the other hand, $$\f(\af_g(a))|_h
= \af_{h\m}(\af_g(a) 1_h) = \af _{h\m}(\af_g (a) 1_g 1_h) = \af
_{h\m}(\af_g (a 1_{g\m} 1_{g\m h})),$$ using again (\ref{triviality}). Now the equality $\bt_g(\f(a))|_h =
\f(\af_g(a))|_h$ follows by (iv) of Definition~\ref{twistdef}, and
(iv) of Definition~\ref{globalization} is proved.

Next we show that

\begin{equation}\label{deltas}
 \f(\D_g) = \f(\A) \cap \bt_g(\f(\A)),
\end{equation}

\noindent for all $g \in G.$ An element from the right hand side
can be written as $\f(a) = \bt_g(\f(b))$ for some $a,b \in \A.$
Then for each $h \in G$ the equality $\f(a)|_h = \bt_g(\f(b))|_h $
means that

\begin{equation*}\label{triviality2}  \af_{h\m}(a 1_h) =
{\widetilde {w}}_{h\m,g} \af_{h\m g}(b 1_{g\m h}){\widetilde
{w}}\m_{h\m,g}.
\end{equation*}

\noindent Taking $h=1$ this gives $a = \af_g(b 1_{g\m}) \in \D_g,$
as ${\widetilde {w}}_{1,g}=1,$  and, consequently, $\f(\D_g) \supseteq
\f(\A) \cap \bt_g(\f(\A)).$ For the reverse inclusion, given an
arbitrary  $a \in \D_g,$ we see using (\ref{triviality}) and (iv)
of Definition~\ref{twistdef} that

\begin{align*}
& \bt _g (\f ({\af}\m _g (a)))|_h =  {\widetilde {w}}_{h\m,g} \,
\af_{h\m g}({\af }\m_{g}(a) 1_{g\m h})\, {\widetilde {w}}\m_{h\m,g} =\\
& \af_{h\m} \circ \af _g ({\af }\m_{g}(a)1_{g\m}  1_{g\m h}) =
\af_{h\m} (a  \af _g (1_{g\m}  1_{g\m h})) = \af_{h\m}(a  1_{ h})=
\f (a) |_h,
\end{align*}
\noindent for all $g, h \in G.$ This yields that $\f (a) =  \bt _g (\f ({\af}\m _g (a)))  $
 for any  $g \in G,$ and $\f(\D_g) \subseteq
\f(\A) \cap \bt_g(\f(\A)).$ Hence (\ref{deltas}) follows and
condition (iii) of Definition~\ref{globalization} is also
satisfied.

 Next we check  that $\f(\A)$ is an ideal in $\B.$
To see this it is enough to show that $\bt_g(\f(a)) \cdot \f(b),
\f(b) \cdot \bt_g(\f(a)) \in \f(\A)$ for all $g \in G$ and $ a,b
\in \A.$  For $h \in G,$   using again (\ref{triviality}) and (iv)
of Definition~\ref{twistdef}, we have

\begin{align*} & \bt_g(\f(a))|_h \cdot \f(b)|_h =
{\widetilde {w}}_{h\m,g}\, \f(a)|_{g\m h} \,  {\widetilde {w}}\m_{h\m,g}
\cdot \f(b)|_h  = \\ & = {\widetilde {w}}_{h\m,g} \, \af_{h\m g} (a
1_{g\m h} ) \,  {\widetilde {w}}\m_{h\m,g} \cdot \af _{h\m}(b 1_h) =
\af_{h\m}(\af_g(a 1_{g\m})1_h) \cdot \af_{h\m}(b 1_h) =
\\& = \af_{h\m} ( \af_g(a 1_{g\m} ) b 1_h) = \f(\af_g(a 1_{g\m}
)b)|_h.
\end{align*}

\noindent  Thus $\bt_g(\f(a)) \cdot \f(b) = \f(\af_g(a 1_{g\m} )b) \in \f(\A)$
and  similarly   $ \f(b) \cdot \bt_g(\f(a))  =\\
\f(b \af_g(a 1_{g\m} )) \in \f(\A),$ as desired.

We show next that our $u$ satisfies (v) of
Definition~\ref{globalization}. For any $g, h, t  \in G,$ by (vi)
of Definition~\ref{twistdef}  we have

\begin{align*} &\f (w_{g,h}) |_t = \af _{t\m} (w_{g,h}1_t) =
w_{t\m,g} \, w_{t\m g,h} \, w\m_{t\m ,gh}  =\\ & 1_{t\m}  1_{t\m
g} 1_{t\m gh} \, \widetilde {w}_{t\m,g} \, \widetilde {w}_{t\m g, h} \,
\widetilde {w}_{t\m,gh} = (\f (1_g 1_{gh}) u_{g,h} ) |_t ,
\end{align*} which gives $$\f (w_{g,h}) =  \f (1_g 1_{gh}) u_{g,h}.$$
The latter equality readily implies $\f (a w_{g,h}) =  \f (a)
u_{g,h}$ for any $a \in {\D}_g {\D}_{gh},$ and the second equality
in (v) of Definition~\ref{globalization} follows similarly.

It remains to prove that
\begin{equation}\label{invariant}
u_{g,h} \, \B = \B = \B \, u_{g,h}
\end{equation} for all $g,h \in G.$

First observe that

\begin{equation*}
u_{g,h} \, \f (a) = \f ({\widetilde{w}}_{g,h} \, a)
\end{equation*} for all $g, h \in G$ and $a \in \A .$ Indeed, for
any $t \in G,$  using (\ref{extend}), one has

\begin{align*} & ( u_{g,h}  \f (a) ) |_t =
{\widetilde{w}}_{t\m,g} \, {\widetilde{w}}_{t\m g,h} \, {\widetilde{w}}\m_{t\m
,gh} \, \af _{t\m}(a 1_t) = \af _{t\m}(w_{g,h} 1_t) \, \af _{t\m}(a 1_t) =
 \f ({\widetilde{w}}_{g,h} \, a) |_t,
\end{align*} as claimed. It follows that
\begin{equation}\label{u.phi}
u_{g,h} \, \f (a)  \in \f (\A ),
\end{equation} for all $g, h \in G.$

It is analogously seen that $u\m_{g,h} \, \f (a) = \f
({\widetilde{w}}\m_{g,h} \,a) ,$ so that
\begin{equation}\label{uinv.phi}
u\m_{g,h} \, \f (a)  \in \f (\A ),
\end{equation} for all $g, h \in G.$

Next we work with the product ${\bt}\m_{t}(u_{g,h}) \, \f (a)$
with arbitrary $g, h, t \in G$ and $a\in \A .$ Note that in view
of (\ref{betacomp}),
$${\bt}\m_{t}(f) = u\m_{t\m,t} \, {\bt}_{t\m} (f) \,  u_{t\m,t}$$ for
any $f \in \F.$ Consequently, applying (\ref{extend}), we have
$${\bt}\m_{t}(u_{g,h}) \, \f (a) = u\m_{t\m,t} \, {\bt}_{t\m}
(u_{g,h}) \,  u_{t\m,t}\, \f (a) = u\m_{t\m,t} \, {u}_{t\m,g} \,
{u}_{t\m g,h} \, {u}\m_{t\m ,gh} \, u_{t\m,t}\, \f (a), $$ which
is contained in $\f (\A )$ by (\ref{u.phi}) and (\ref{uinv.phi}).
This yields that ${\bt}\m_{t}(u_{g,h}) \, \f ( \A ) \subseteq \f
(\A ) ,$ and applying $\bt _ t$ we obtain that $$u_{g,h} \,
{\bt}_{t}(\f ( \A )) \subseteq {\bt}_{t}( \f (\A )) ,$$ for all
$g,h, t \in G.$ Hence   $u_{g,h} \, \B \subseteq  \B ,$ and one
similarly shows that $u\m_{g,h} \, \B \subseteq  \B .$ An
analogous argument gives  $ \B \, u_{g,h}, \,  \B \, u\m_{g,h}
\subseteq \B ,$  and (\ref{invariant}) follows. \fim \\

\end{section}

\begin{section}{Some remarks on transitive twisted partial actions}\label{sec:transitive}

In this section we start our treatment of twisted partial actions
on products of indecomposable rings by making some preliminary
remarks. Let
\begin{equation}\label{decomp1} \A = \prod _{\lb \in
\Lb} {\R}_{\lb},
\end{equation} where each ${\R}_{\lb} $ is an
indecomposable unital ring and $\Lb $ is some
 non-necessarily finite index set. It is directly verified that a decomposition of $\A $
 into a product of indecomposable factors (blocks of $\A $ ) is unique up to a permutation of indices.
 Clearly, each ideal of $\A $ which is a direct factor, i.e. is
 generated by an idempotent which is central in $\A ,$ has to be a
 product of some blocks.  Thus, if a twisted partial action $\af $ of a group $G$ on $\A $
 is given such that each $\D _g$ is a unital ring, then for every $g \in G$ the ideal $\D _g$ is a
 product of some ${\R }_{\lb}$'s. It follows that $\af $ permutes the
  indecomposable factors in the sense that  if  $g \in G$ and  $\lb \in \Lb $ are such that
  $\R _{\lb } \subseteq \D _{g\m},$ then $\af _g (\R _{\lb}) = \R_{{\lb}'}$ for some ${\lb}' \in \Lb .$ Note
  also that since each ${\D}_g $ is unital, $w_{g,h} \in \U
({\D}_g {\D}_{gh})$ for any $g,h \in G,$ and hence $w_{g,h}
{\R}_{\lb }$ is either $0$ or $  {\R}_{\lb }$ for any $\lb \in \Lb
,$ and similarly for $ {\R}_{\lb }w_{g,h}.$ We shall say  that
$\af $ is {\it transitive} if for any
  $\lb , {\lb}' \in \Lb $ there exists  $g \in G$ such that  $\R _{{\lb}} \subseteq \D _{g\m}$ and
  $ \af _g (\R _{\lb} ) = \R _{{\lb}'}.$ Fix an arbitrary block  from (\ref{decomp1}) which
  shall be denoted by  ${\R }_1.$ It is readily verified that
 $\af $ is  transitive if and only if for any
  $\lb \in \Lb $ there exists  $g \in G$ such that  $\R _1 \subseteq \D _{g\m}$ and
  $ \af _g (\R _1 ) = \R _{\lb}.$ \\

  Suppose that $\af $ is transitive and set
  $$H = {\rm St}_G (\R _1 )= \{g \in G : \R _1 \subseteq \D _{g\m}, \af _g (\R _1 ) = \R _1
  \},$$ the stabilizer of $\R _1 $ in $G ,$ and let ${\Lb}' \ni 1$ be a
  left transversal of $H$ in $G,$ i.e. $G = \bigcup _{g \in
  {\Lb}'} gH,$ a disjoint union. One may evidently assume that $\Lb
  \ni 1$ is a subset of ${\Lb}'$ so that  $\R _1 \subseteq \D _{g\m}$ for all $g \in \Lb $ and
  (\ref{decomp1}) can be
  rewritten as

\begin{equation}\label{decomp2} \A = \prod _{g \in
\Lb} {\R}_{g},
\end{equation} where ${\R}_g = \af _g (\R _1).$\\

For $x \in G$ denote by $\bar x$ the element of ${\Lb }'$ with
$\bar x H = x H.$ Keeping our notation we shall use the next easy
properties.

\begin{lema}\label{easy} For $g, x \in G$  we have:\\

\noindent (i) $g\in {\Lb }', \R _1 \subseteq \D _{g\m}  \iff  g\in \Lb  ;$\\

\noindent  (ii) $g , \, {\overline{xg}} \in \Lb \iff {\R}_g
\subseteq {\D}_{x\m},$ and if this holds then
$\af _x ({\R}_{g})= {\R}_{\overline{xg}};$\\

\noindent  (iii) $ g\in \Lb ',  {\R}_{\overline{x\m g}} \subseteq {\D}_{x\m}
\Longrightarrow g\in \Lb .$
\end{lema}

\p (i) We only need to see the ``$\Longrightarrow $'' part, so let $ \R _1 \subseteq \D _{g\m}$
for some
$g\in {\Lb }'.$ Then $\af _g (\R _1)$ must be a block of $\A ,$ so it equals $\R _t $ for some
$t \in \Lb .$
Hence
\begin{align*}
& \R _1 = {\af}\m _g  \circ \,  \af _t (\R _1) = w\m_{g\m ,g } \, {\af }_{g\m } \circ {\af }_t (\R _1)  \,
w_{g\m ,g}
=\\ &  w\m_{g\m ,g } \,w_{g\m , t} \, \af _{g\m t}( {\R}_1) \,  w\m_{g\m ,t} \, w_{g\m ,g} =
 \af _{g\m t}( {\R}_1) .
\end{align*}  Consequently, $g\m t \in H ,$ and $g =t  \in \Lb .$\\

(ii) Let $g ,  {\overline{xg}} \in \Lb .$ Then $\R _1 \subseteq
\D _{(\overline{xg})\m} \cap \D _{g \m},$ and since
$(\overline{xg})\m xg \in H,$ one has that
${\af}_{\overline{xg}} \circ {\af}_{ (\overline{xg})\m xg}$ is
applicable to ${\R}_1$ and by Definition~\ref{twistdef} so too is
 ${\af}_{xg} = {\af}_{\overline {xg} \cdot (\overline{xg})\m xg}.$
Thus $\R _1 \subseteq \D _{({xg})\m}$ and  using again
Definition~\ref{twistdef} we see that $\R _g = \af _g ({\R}_1)
\subseteq \af _g (\D _{(xg)\m} \cap \D _{g\m }) = \D _{x\m} \cap
\D _{g },$ so that $\af _x$ is applicable to ${\R}_g.$  Moreover,
 \begin{align*}
 & \af _{x} ({\R}_g) = w_{x,g} \af _{xg} ({\R}_1) w\m_{x,g} =
  w_{x,g}\,   w\m_{\overline{xg}, (\overline{xg})\m xg}\,
 {\af}_{\overline {xg}} \circ {\af }_{(\overline{xg})\m xg} ({\R}_1)\,
 w_{\overline{xg}, (\overline{xg})\m xg}\, w\m_{x,g}\\ &
= w_{x,g}\, w\m_{\overline{xg}, (\overline{xg})\m xg}\,
{\af}_{\overline {xg}} ({\R}_{1}) \, w_{\overline{xg},
(\overline{xg})\m xg}\, w\m_{x,g} \subseteq {\R}_{\overline {xg}},
\end{align*} and consequently $\af _{x} ({\R}_g) = {\R}_{\overline
{xg}},$ as $\af _{x} ({\R}_g)$ must be a block.

Conversely,
\begin{align*}
& {\R}_g \subseteq {\D}_{x\m} \Longrightarrow g\in \Lb, {\R}_g
\subseteq {\D}_{x\m} \cap {\D}_{g} = {\D}_{x\m} {\D}_{g}  \Longrightarrow  \\
& {\R}_1 =
{\af}\m_g ({\R}_g) \subseteq  {\af}\m_{g} ({\D}_{x\m} {\D}_{g} )  \subseteq {\D}_{g\m x\m}
\Longrightarrow
  {\R}_1  \subseteq  {\D}_{g\m x\m}  \cap {\D}_{g\m x\m \overline{xg}} \Longrightarrow \\
& {\R}_1 = \af _{ (\overline{xg})\m xg } \, (\R _1 ) \subseteq   {\af}_{(\overline {xg} )\m xg}
({\D} _{g\m x\m}  \cap {\D}_{g\m x\m \overline {xg} })
\subseteq {\D}_{(\overline {xg} )\m},
\end{align*} which gives $\overline {xg} \in \Lb .$\\

(iii) If $\R _ {\overline{x\m g}} \subseteq {\D}_{x\m }$ then by (ii) we have that
$\af _x (\R _ {\overline{x\m g}} ) =
 (\R _{\overline{x x\m g}}) =  {\R}_g ,$ i.e. $g \in \Lb .$ \fim \\

 Note that taking $g=1$ in (ii) of  Lemma~\ref{easy}, we have
 that

\begin{equation}\label{fromeasy}
\R _1 \subseteq \D _{x\m } \iff \R _1 \subseteq \D _{{\bar x}\m },
 \end{equation} for any $x \in G.$\\

  Given $g \in \Lb ,$  we denote by ${\rm pr}_g : \A \to \R _g  $ the projection map; the element
  ${\rm pr}_g (a)$ shall be called the $g$-{\it entry} of $a.$
  Using this notation one may write $$ a = \prod_{g\in \Lb } {\rm pr}_g (a)$$ for any $a \in \A .$
  We shall work with multiplicative maps   ${\0}_x :
\A \to \A ,$ $x \in G,$ defined as follows:
for arbitrary $a \in \A$ let   ${\0}_x (a) = 1_{\A }$  if $\overline {x\m} \notin \Lb ,$
and if $\overline {x\m} \in \Lb ,$ 3 $\R _1 \subseteq {\D}_x,$ then
$${\rm pr}_g \,  {\0}_x (a) = \left\{\begin{array}{cl} 1_{\R _g }, \,
& \mbox{if} \, \, g \neq \overline {x\m}, \\ {\rm pr}_{g} \, {\af
}\m  _x ({\rm pr}_1 a), & \mbox{otherwise.} \end{array} \right. $$
\\ We shall frequently use the next obvious equality:
\begin{equation}\label{obvious}    {{\0} _x} (a) =   {{\0} _x} (a 1_h)
\end{equation} for any $x \in G,$ $a \in \A $ and $h \in H= {\rm St}_G (\R _1 ).$ We also need the
following properties.

\begin{lema}\label{properties}  For any $x \in G$ and $g \in {\Lb }'$ we have\\

\noindent  (i) $1_x \, {{\0}_{g\m}} \,(a 1_{g\m x}) =  1_x \, {\0}_{g\m} \, (a ),$ for all $a \in \A ;$\\

\noindent (ii) $1_x    \,  {\0}_{g\m }    \circ \af _{g\m} \circ  \af _x \circ
{\af}\m_{({\overline{x\m g}})\m} (a) =
 \af _x (1_{x\m } \,  {\0}_{({\overline{x\m g}})\m} (a) ),$
 for all \\ $  a\in {\D }_{({\overline{x\m g}})\m x\m} \cap {\D }_{({\overline{x\m g}})\m} \cap
 {\D }_{({\overline{x\m g}})\m x\m g}.$
\end{lema}

\p (i) Since $\0 _{g\m} $ is multiplicative, it is enough to check
that $1_ x \0 _ {g\m} (1_{g\m x}) = 1_x.$ If $g \notin \Lb ,$ then
$\0 _{g\m} (1_{g\m x}) = 1_{\A},$ and our equality is clearly
satisfied. Moreover, the equality $\0 _{g\m } (1_{g\m x}) =
1_{\A}$ also holds if $g \in \Lb$ and $\R _1 \subseteq \D _{g\m
x}.$ If $g \in \Lb$ and $\R _1 \not \subseteq   \D _{g\m x},$ then
the $g$-entry of $\0 _{g\m } (1_{g\m x})$ is $0$ and all other
entries are $1$'s. By (\ref{fromeasy}) and Lemma~\ref{easy}, $\R
_g \not \subseteq \D _x ,$ and consequently $1_x \0 _{g\m } (1_{g\m x}) = 1_x $ holds also in this case.\\

(ii) Note that ${\D }_{({\overline{x\m g}})\m x\m} \cap {\D
}_{({\overline{x\m g}})\m} \cap
 {\D }_{({\overline{x\m g}})\m x\m g}$ is the set of all $a\in \A $ for which
  $ \af _{g\m} \circ  \af _x \circ  {\af}\m_{({\overline{x\m g}})\m}$ is applicable. Denote by
  $a_1$ the left hand side of the equality in (ii) and by $a_2$ its right hand side.
It is easy to see that if $g \notin \Lb $ or  $\overline{x\m g}
\notin \Lb ,$ then $a_1 = 1_x  = a_2,$ so suppose  that $g,
\overline{x\m g} \in \Lb .$ Evidently $a_1, a_2 \in \D _x$ and by
the definition of the $\0 $'s, for any $t \in \Lb $ with $t \neq
g$ and $ \R _t \subseteq \D _x ,$ we have ${\rm pr}_t (a_1) =
1_{\R _t} =  {\rm pr}_t (a_2) .$ Furthermore,  by
Lemma~\ref{easy}, $\R _g \subseteq \D _x$ and 
\begin{align*}& {\rm pr}_g (a_1) =
{\rm pr}_g( {\0}_{g\m }    \circ \af _{g\m} \circ  \af _x \circ
{\af}\m_{({\overline{x\m g}})\m}  (a) )=  {\rm pr}_g(    \af _x
\circ  {\af}\m_{({\overline{x\m g}})\m}  (a) ) =\\& {\rm pr}_g( \af
_x (1_{x\m } \,  {\0}_{({\overline{x\m g}})\m} (a) )) =  {\rm
pr}_g (a_2),
\end{align*}  and consequently $a_1 = a_2.$ \fim

\end{section}

\begin{section}{Going to an equivalent twisted partial action}\label{sec:goingto}

In this section we adapt to the partial and non-commutative
context the idea of the corestriction from Homological Algebra in
order to produce a more manageable  twisted partial action,
equivalent to the initial one.  We give the following:

\begin{defi}\label{equivalent} Two twisted partial actions of $G$
on a ring $\A $
$${\af } = (\{{\D }_x \}_ {x \in G}, \{{\af}_x \}_ {x \in G},
 \{w_{x,y} \}_ {(x,y) \in G\times G})$$ and
 $${\af }' = (\{{\D}_x \}_ {x \in G}, \{{\af}'_x \}_ {x \in G},
 \{w'_{x,y} \}_ {(x,y) \in G\times G}),$$ such that each $\D _x$ $(x\in G)$ is a unital ring with unity
 $1_x ,$ shall be called
 equivalent if there exists a function $$ G \ni x \mt {\e}_x \in
 {\U }(\D _x) \subseteq \A $$ such  that for all $x, y \in G$ and
 $a \in \D _{x\m} $ we have
\begin{equation}\label{alpha}
{\af}'_x (a) = {\e}_x {\af}_x (a) {\e}\m
 _x,
\end{equation} and
 \begin{equation}\label{w} w'_{x,y} = \e _x \,  {\af}_x (  \e _y 1_{x\m}) \, w_{x,y} \, {\e}\m_{xy}.
\end{equation}
\end{defi}

Observe that in the above definition the twisted partial actions
have common domains $\D _x.$ It is readily checked that the
definition is correct, i.e. our relation is reflective, symmetric
and transitive.\\

\begin{lema}\label{epsilon} Let $${\af } = (\{{\D }_x \}_ {x \in G}, \{{\af}_x \}_ {x \in G},
 \{w_{x,y} \}_ {(x,y) \in G\times G})$$ be a twisted partial action of $G$ on $\A $ such that each
 $\D _x$ $(x\in G)$ is a unital ring with unity
 $1_x .$ Given a function  $$ G \ni x \mt {\e}_x \in
 {\U }(\D _x) \subseteq \A ,$$ define ${\af}'_x$ $(x \in G)$ by (\ref{alpha}) and $w_{x,y}$ $(x,y \in G)$
 by (\ref{w}). Then
$${\af }' = (\{{\D}_x \}_ {x \in G}, \{{\af}'_x \}_ {x \in G},
 \{w'_{x,y} \}_ {(x,y) \in G\times G})$$ is a twisted partial action of $G$ on $\A .$
\end{lema}

\p The items (i), (ii), (iii) and (v) of Definition~\ref{twistdef} are immediate, and the
verification of (iv) is straightforward.
As to (vi), for arbitrary $x, y, z \in G$ we have
\begin{align*}  & {\af}'_x(1_{x\m } w'_{y,z})w'_{x,yz} =
\e _x \af _x (1_{x\m} \e _y \af _y (\e _z 1_{y\m} ) w_{y,z} {\e}\m_{yz} ) {\e}\m _x
\e _x \af _x ( \e _{yz} 1_{x\m }) w_{x,yz} {\e}\m_{xyz} \\
& = \e _x  \af _x ( 1_{x\m } \e _y) ({\af}_x \circ {\af}_y (\e _z 1_{y\m x\m} 1_{y\m}) )
w_{x,y} w_{xy,z} {\e}\m_{xyz},
\end{align*} by cancelling the $\e $'s, and using  the $2$-cocycle equality for the $w$'s and
(\ref{triviality}). Now composing ${\af}_x \circ {\af}_y$ according to
(iv) of Definition~\ref{twistdef}, we easily come to $1_x w'_{x,y} w'_{xy,z} ,$ as claimed.  \fim \\

In what follows in this section we assume all notation from the
previous one, in particular, $\A $ is a product of blocks as in
(\ref{decomp2}), and $\af $ is a transitive twisted partial action
of $G$ on $\A.$ We remind that we shall  write the element
$w_{x,y}$ also as $w[x,y],$ the standard notation used in
Homological Algebra for $2$-cocycles. Since each $\D _x$ has unity
$1_x ,$ we can rewrite the
$2$-cocycle equality (vi) from Definition~\ref{twistdef} as\\

\begin{equation}\label{cocycle}
\af_x(1_{x\m} w[y,z]) \, w[x, yz] =  w[x,y] \,  w[xy,z], \;
\forall x,y,z \in G.
\end{equation}\\

Set

\begin{equation}\label{corestr}
w'[x,y ] = 1_x 1_{xy} \,  \prod_{g \in \Lb} {\0 }_{g\m }(  w[g\m x
\cdot \overline{x\m g}, \, (\overline{x\m g})\m \cdot y \cdot
\overline{y\m x\m g}]).
\end{equation}

Since $\R _1 \subseteq \D _h$ for any $h \in H = {\rm St}_G (\R
_1),$ and $$w[g\m x \cdot \overline{x\m g}, \, (\overline{x\m
g})\m \cdot y \cdot \overline{y\m x\m g}] \in {\U}({\D}_{g\m x
\cdot \overline{x\m g}} \, {\D}_{g\m x \cdot y \cdot
(\overline{y\m x\m g})\m} )$$ with $g\m x \cdot \overline{x\m g},
g\m x \cdot y \cdot (\overline{y\m x\m g})\m \in H,$ it follows
that $w'[x,y ]$ is an invertible element in  $\D _x \D _{xy}$ for
all $x,y \in G.$\\

Our goal in this section is the next:

\begin{prop}\label{corestr-equiv} With the above notation,
$${\af }' = (\{{\D }_x \}_ {x \in G}, \{{\af}'_x \}_ {x \in G},
 \{w'_{x,y} \}_ {(x,y) \in G\times G}) $$ is a transitive twisted
 partial action of $G$ on $\A $ with
 \begin{equation}\label{equiv}
 w[x,y]= {\af }_x (\e _y 1_{x\m}) \e _x \, w'[x,y] \, {\e }\m_{xy}, \, \, \,
\forall x,y \in G,
\end{equation} where ${\af }'_x (a) = {\e}\m_x {\af }_x (a) {\e}
_x$   and
$$ \e _x = 1_x \prod _{g \in \Lb} {\0}_{g\m} (w[g\m, x] \, w[g\m x
\cdot \overline{x\m g}, \, (\overline{x\m g})\m ]\m) \in {\U}(\D
_x),$$ for all $x \in G$ and  $a \in \D _{x\m}.$ In particular,
$\af $ and ${\af }'$ are equivalent.
\end{prop}

\p First note that by  Lemma~\ref{easy} and (\ref{fromeasy}), $$
\R _g \subseteq \D _x \iff \R _1 \subseteq \D _{g\m x},$$ and
since $$w[g\m, x] \in {U}(\D _{g\m} \D _ {g\m x}),  \,  w[g\m x
\cdot \overline{x\m g}, \, (\overline{x\m g})\m ] \in {U}(\D _{g\m
x \cdot \overline{x\m g} } \, \, \D _ {g\m x})$$ and $g\m x \cdot
\overline{x\m g} \in H,$ it follows that $\e _x $ is an invertible
element in $\D _x .$ Thus by Lemma~\ref{epsilon}, ${\af }'$ is a
twisted partial action, and it is clearly transitive, so it
remains to concentrate on proving (\ref{equiv}). We have

\begin{align*} & w[x,y] = 1_x 1_{x y} \prod _{g\in \Lb} {\0}_{g\m}
({\af}_{g\m}(w[x,y]1_g)),
\end{align*} for any $x, y \in G$ by the definition of ${\0}_{g\m}.$ Applying the
$2$-cocycle equality (\ref{cocycle}) for the triple $(g\m, x,y),$
and using (\ref{obvious}), we obtain
\begin{align*} & w[x,y] = 1_x 1_{xy} \prod _{g\in \Lb} {\0}_{g\m }
(w[g\m ,x]w[g\m x, y]\, 1_{g\m x \cdot \overline{x\m g}} \, 1_{g\m
x \cdot y \cdot \overline{y\m x\m g}} \,  \, w[g\m ,xy]\m ),
\end{align*} as $ g\m x \cdot \overline{x\m g}, g\m
x \cdot y \cdot \overline{y\m x\m g} \in H. $ Next, equality
(\ref{cocycle}) used for  the triples $$(g\m x \cdot \overline{x\m
g}, \, (\overline{x\m g})\m, y)\, \, \,  \mbox{and} \, \, \, (g\m
x \cdot \overline{x\m g}, \, (\overline{x\m g})\m \cdot y \cdot
\overline{y\m x\m g}, \, (\overline{y\m x\m g})\m  )$$
respectively gives
\begin{align*} & w[g\m x, y]\,
1_{g\m x \cdot \overline{x\m g}}  =  w[g\m x \cdot \overline{x\m
g}, (\overline{x\m g})\m ]\m \cdot \\ & \af _{g\m x \cdot
\overline{x\m g}} \,(w[ (\overline{x\m g})\m , y] 1_{(
\overline{x\m g})\m \cdot x\m g}) \,w[g\m x\cdot \overline{x\m g},
( \overline{x\m g})\m \cdot y], \end{align*} and
\begin{align*} & w[g\m x \cdot \overline{x\m g}, (\overline{x\m g})\m \cdot
y] \, 1_{g\m x y \cdot \, \overline{y\m x\m g}} =  \\ & \af _ {
g\m x \cdot \,  \overline{x\m g} } \, ( w[(\overline{x\m g})\m
\cdot y \cdot \overline{y\m x\m g}, \, (\overline{y\m x\m g})\m ]
\, 1_{(\overline{x\m g})\m \cdot \, x\m g})\m \\ &  \cdot w[g\m x
\cdot \overline{x\m g} , \, (\overline{x\m g})\m \cdot y \cdot
\overline{y\m x\m g} ] \,  w[g\m x y \cdot \overline{y\m x\m g},
(\overline{y\m x\m g})\m].
\end{align*} Making consecutive replacements in our expression of
$w[x,y],$ we come to
\begin{align*}\label{a1} w[x,y] =\e _x
\, (\prod _{g\in \Lb} {\0}_{g\m} \circ \af _{g\m x \cdot
\overline{x\m g}} \, ( a  \, 1_{(\overline{x\m g})\m \cdot x\m g}
) )\,  w'[x,y] {\e }_{xy},
\end{align*} where
\begin{align*} & a =  \af _{g\m x \cdot \, \overline{x\m g}} \, ( w[(\overline{x\m
g})\m , y]  , w[(\overline{x\m g})\m \cdot y \cdot \overline{y\m
x\m g}, \, (\overline{y\m x\m g})\m  ]\m \\ & 1_{(\overline{x\m
g})\m \cdot x\m g}) \in {\D}_{(\overline{x\m g})\m } \,
{\D}_{(\overline{x\m g})\m \cdot y} \, {\D}_{(\overline{x\m g})\m
\, \cdot y \cdot  \,  \overline{y\m x\m g}},
\end{align*} so that it remains to show that

\begin{equation}\label{a2}
1_{xy} \, {\e }_x \prod _{g\in \Lb } {\0}_{g\m } \circ \af _{g\m x
\cdot \overline{x\m g}} \, ( a  \, 1_{(\overline{x\m g})\m \cdot
x\m g} )  = {\af }_{x} ( {\e}_y 1_{x\m }) \e _x,
\end{equation} for any $x, y \in G.$

It will be convenient to take in the above equality $g$ running
over ${\Lb}';$ it does not cause any harm by the definition of
${\0}_{g\m}.$  Since $a$ lies in $ {\D}_{(\overline{x\m g})\m}$ we
may write $$a  \, 1_{(\overline{x\m g})\m \cdot x\m g} = \af
_{(\overline{x\m g})\m} \circ {\af}\m_{(\overline{x\m g})\m} \, (
a \, 1_{(\overline{x\m g})\m \cdot x\m g}),$$ so that when using
this to the left hand side of (\ref{a2}) and composing $\af _{g\m
x \cdot \overline{x\m g}} \circ \af _{(\overline{x\m g})\m},$ we
obtain that
\begin{align*} & 1_{xy} \, {\e }_x \prod _{g\in {\Lb }' } {\0}_{g\m } \circ \af _{g\m x \cdot
\overline{x\m g}} \, ( a  \, 1_{(\overline{x\m g})\m \cdot x\m g}
) =  1_x 1_{xy}  \prod _{g \in {\Lb }'} {\0}_{g\m} (w[g\m, x] \,
\af _{g\m x } \circ \\ & {\af }\m_{ (\overline{x\m g})\m} \, ( a
\, 1_{(\overline{x\m g})\m \cdot x\m g} ) \,   w[g\m x \cdot
\overline{x\m g}, \, (\overline{x\m g})\m ]\m),
\end{align*} since ${\af }\m_{ (\overline{x\m g})\m} \, ( a \,
1_{(\overline{x\m g})\m \cdot x\m g} ) \in {\D}_{\overline{x\m g}
} \, {\D}_{g\m x}.$  Next, putting $x=1$ in (i) of
Lemma~\ref{properties} and using (\ref{triviality}), we have that
the left hand side of (\ref{a2}) equals
\begin{align*}  & 1_x 1_{xy} \prod _{g \in {\Lb }'} {\0}_{g\m} (w[g\m, x] \, \af _{g\m
x } (1_{x\m}  {\af }\m_{ (\overline{x\m g})\m} \, ( a \,
1_{(\overline{x\m g})\m \cdot x\m g} )) \\ &  w[g\m x \cdot
\overline{x\m g}, \, (\overline{x\m g})\m ]\m) = \\ &  1_{xy}
\prod _{g \in {\Lb }'} {\0}_{g\m} (  \af _{g\m }\circ {\af }_{x }
(1_{x\m} {\af }\m_{ (\overline{x\m g})\m} \, ( a \,
1_{(\overline{x\m g})\m \cdot x\m g} )))   {\e }_x = \\ & 1_{xy}
\prod _{g \in {\Lb }'} {\0}_{g\m} (  \af _{g\m }\circ {\af }_{x }
\circ  {\af }\m_{ (\overline{x\m g})\m} \, ( a \,
1_{(\overline{x\m g})\m \cdot x\m } 1_{(\overline{x\m g})\m \cdot
x\m g} )) {\e }_x.
\end{align*} The latter, using both items of
Lemma~\ref{properties}, equals
\begin{align*} & 1_{xy} \prod _{g \in
{\Lb }'}  {\af }_{x } (1_{x\m }  {\0 }\m_{ (\overline{x\m g})\m}
\, ( a \, 1_{(\overline{x\m g})\m \cdot x\m } 1_{(\overline{x\m
g})\m \cdot x\m g} )) {\e }_x= \\ &  {\af }_{x } (1_{x\m } 1_{y}
\prod _{g \in {\Lb }'}  \, {\0 }\m_{ (\overline{x\m g})\m} \, ( a
)) {\e }_x= \af _x (\e _y 1_{x\m }) \e _x,
\end{align*} as $\overline{x\m g}$ runs over ${\Lb }'.$ \fim

\end{section}

\begin{section}{Existence of a globalization}\label{sec:existence}

In this section we prove the  existence of globalizations of
partial actions on rings which are products of blocks. According
to Theorem~\ref{globreduced}, it is enough to construct $\widetilde
{w}[x,y] = \widetilde {w}_{x,y} \in {\U}(\A)$ which satisfy the
extended partial $2$-cocycle equality (\ref{extend}). For this
purpose, keeping the notation of the previous sections, define for
every $x \in G$ the map $\tilde{\af }_x : \A \to \A ,$ by
\begin{equation*}
\tilde{\af }_{x}(a) =   \af _x (a 1_{x\m }) + 1_{\A }- 1_{x},
\end{equation*} with  $a \in \A ,$ and  set
\begin{equation}\label{tildew}
\widetilde{w }[x,y] =  \tilde{\af }_x (\tilde{\e}_y) \tilde{\e }_x
\widetilde{w'}[x,y ] \, \tilde{\e }_{xy}\m \in \U (\A ),
\end{equation} where

\begin{equation*} \tilde{\e}_x = {\e}_x + 1_{\A } - 1_x \in \U (\A ),
\end{equation*} and
\begin{equation*}\label{tildecorestr}
\widetilde{w'}[x,y ] =   \prod_{g \in {\Lb}'} {\0 }_{g\m }( w[g\m x
\cdot \overline{x\m g}, \, (\overline{x\m g})\m \cdot y \cdot
\overline{y\m x\m g}]) \in \U (\A ),
\end{equation*}  with $x, y \in G.$ We have  by
(\ref{triviality}) that $1_x 1_{xy} \tilde{\af }_x (\tilde{\e}_y)
= {\af }_x (1_{x\m }1_{y}\tilde{\e}_y)  = {\af }_x (1_{x\m
}{\e}_y)$ and consequently
\begin{equation*} 1_x 1_{xy} \widetilde{w }[x,y] = w[x,y],
\end{equation*} for all $x,y \in G .$\\

We shall use the following composition law and  $2$-cocycle
equality:
\begin{lema}\label{twolaws} For every $x, y, z \in G$ and $a \in \A$ we have
\begin{equation}\label{composition} {\af}_x (1_{x\m } \tilde{\af }_y (a)) = 1_x
\widetilde{w}[x,y] \, \tilde{\af }_{xy} (a) \,  \widetilde {w} [x,y]\m .
\end{equation} and
\begin{equation}\label{cocycle2} {\e }\m_x \af _x (1_{x\m}
\widetilde{w'}[y,z] ) \e _x \, \widetilde{w'}[x,yz] = 1_x \widetilde{w'}[x,y]
\, \widetilde{w'}[xy,z].
\end{equation}
\end{lema}

\p Notice first that multiplying the left hand side of
(\ref{composition})  by $1_{xy}$ gives
\begin{align*} & {\af}_x (1_{x\m } 1_y \tilde{\af }_y (a )) =
{\af}_x ( {\af }_y (a 1_{y\m x\m} 1_{y\m} )) = w[x,y] \m
{\af}_{xy} ( a 1_{y\m x\m} ) \m w[x,y]\m = \\ & 1_x 1_{xy}
\widetilde{w}[x,y] \, \tilde{\af }_{xy} (a) \,  \widetilde {w} [x,y]\m ,
\end{align*} which is the right hand side  multiplied by $1_{xy}.$ On the
other hand, it is directly checked that both sides  multiplied by
$1_{\A } - 1_{xy}$ equal $1_x (1_{\A } - 1_{xy}),$ which proves
(\ref{composition}).\\

To prove (\ref{cocycle2}),  we first observe   using (\ref{obvious}),
(i) of Lemma~\ref{properties} and the obvious equality $\overline
{t_1 \overline{t_2}} = \overline{t_1 t_2},$ $(t_1, t_2 \in G),$
that
\begin{align*} &  \af _x  (1_{x\m } \widetilde{w'}[y,z]) = \prod _{g \in {{\Lb}'}} \af _x (1_{x\m} \,
 \0 _{(\overline{x\m
g})\m} (1_{(\overline{x\m g})\m \cdot x\m}1_{(\overline{x\m g})\m}
1_{(\overline{x\m g})\m x\m g} \, a)),
\end{align*} where $$a = w[(\overline{x\m
g})\m \cdot y \cdot \overline{y\m x\m g}, (\overline{y\m x\m g})\m
\cdot z \cdot \overline{z\m y\m x\m g} ],$$ taking into account
that   $\overline{ x\m g}$ runs over ${\Lb}'$ when so too does $g$
with fixed $x \in G .$  Next, by (ii) of Lemma~\ref{properties},
(\ref{triviality}) and the composition law (iv) of
Definition~\ref{twistdef}, we obtain
\begin{align*} & \tilde{\e }\m _x \af _x  (1_{x\m } \widetilde{w'}[y,z]) \tilde{\e
}_x=\\ & \tilde{\e }\m _x ( \prod _{g \in {{\Lb}'}} 1_x  \, \0
_{g\m} \circ \af _{g\m} \circ {\af}_x \circ
 {\af }\m _{(\overline{x\m
g})\m} (1_{(\overline{x\m g})\m \cdot x\m} 1_{(\overline{x\m
g})\m} 1_{(\overline{x\m g})\m x\m g} \, a) )\tilde{\e }_x= \\
&  1_x \prod _{g \in {{\Lb}'}} \0 _{g\m} ( w[g\m x \cdot
\overline{x\m g}, \, (\overline{x\m g})\m ]   \af _{g\m x} \circ
\\ &
 {\af }\m _{(\overline{x\m
g})\m} (1_{(\overline{x\m g})\m \cdot x\m} 1_{(\overline{x\m
g})\m} 1_{(\overline{x\m g})\m x\m g} \, a) \, w[g\m x
\cdot \overline{x\m g}, \, (\overline{x\m g})\m ]\m ) = \\
&  1_x \prod _{g \in {{\Lb}'}} \0 _{g\m} ( w[g\m x \cdot
\overline{x\m g}, \, (\overline{x\m g})\m ] \,    \af _{g\m x
\cdot \overline{x\m g} \cdot  (\overline{x\m g})\m} \circ
\\ &
 {\af }\m _{(\overline{x\m
g})\m} ( 1_{(\overline{x\m g})\m} 1_{(\overline{x\m g})\m x\m g}
\, a) \, w[g\m x
\cdot \overline{x\m g}, \, (\overline{x\m g})\m ]\m ) = \\
&  1_x \prod _{g \in {{\Lb}'}} \0 _{g\m} (  \af _{g\m x \cdot
\overline{x\m g} } \, (  1_{(\overline{x\m g})\m x\m g} \, a)).
\end{align*} Now,  ${\0}_{g\m}(b) = {\0}_{g\m}({\rm pr}_1 b)$ $(b \in \A ),$
by the definition of ${\0 }_{g\m} ,$ and therefore,
\begin{align*} & {\e }\m_x \af _x (1_{x\m}
\widetilde{w'}[y,z] ) \e _x \, \widetilde{w'}[x,y z] = \\ & 1_x \prod _{g
\in {\Lb}'} \0 _{g\m} ({\rm pr}_1 (  \af _{g\m x \cdot
\overline{x\m g} } \,  ( {\rm pr}_1 a )   w[g\m x \cdot
\overline{x\m g}, (\overline{x\m g})\m \cdot y z \cdot
(\overline{z\m y\m x\m g})\m])),
\end{align*} because $g\m x \cdot \overline{x\m g} \in H$ and ${\R
_1}\subseteq {\D}_{(\overline{x\m g})\m \cdot x\m g}.$ Finally,
the automorphisms ${\rm pr}_1 \circ {\af }_h \circ {\rm pr}_1, $
$h \in H$ of ${\R }_1$ and the twisting  ${\rm pr}_1 (w[h, h']),$
$h \in H$ form a twisted global action of $H$ on ${\R }_1,$ and
using the (global) $2$-cocycle equality for the triple
$$(g\m x \cdot \overline{x\m g}, (\overline{x\m g})\m \cdot y
\cdot \overline{y\m x\m g} , (\overline{y\m x\m g})\m \cdot z
\cdot \overline{z\m y\m x\m g} ),$$ in which each entry belongs to
$H,$ we conclude that

\begin{align*} & {\e }\m_x \af _x (1_{x\m}
\widetilde{w'}[y,z] ) \e _x \, \widetilde{w'}[x,y z] = \\ &  1_x \prod _{g
\in {\Lb}'} \0 _{g\m} ({\rm pr}_1 (  w[g\m x \cdot \overline{x\m
g}  , (\overline{x\m g})\m \cdot y \cdot \overline{y\m x\m g}] \\
&  \cdot w[g\m x\cdot y \cdot \overline{y\m x\m g}, (\overline{y\m
x\m g})\m \cdot z \cdot \overline{z\m y\m x\m g} ]) ) = 1_x
\widetilde{w'}[x,y] \, \widetilde{w'}[xy,z],
\end{align*} as desired. \fim \\

Now we state the following:

\begin{teo}\label{exist} Let $\A $ be a unital ring  which is a (non-necessarily finite)
product of indecomposable rings. A twisted partial action
$$\af = (\{\D_x \}_ {x \in G}, \{{\af}_x \}_ {x \in G},
 \{w[x,y] \}_ {(x,y) \in G\times G}),$$ of a group
$G$ on  $\A $  is globalizable exactly when each ${\D}_x$ $(x \in
G)$ is a unital ring.
\end{teo}

\p It was observed already that the existence of a globalization
implies that each $\D _x, x \in G$ is unital. For the converse let
$1_x $ be the unity of $\D _x $ and assume first that $\af $ is
transitive. Keeping the above notation we show that the $\widetilde
{w}[x,y]$'s defined in (\ref{tildew}) satisfy
\begin{equation}\label{extend2} \af_x( \widetilde{w}[y,z]\, 1_{x\m})\, \widetilde{w}[x,
yz] =1_x \, \widetilde{w}[x,y] \, \widetilde{w}[xy,z],
\end{equation} for any $x, y, z \in G.$

The left hand side of (\ref{extend2}) equals
\begin{align*}  & \af _x (1_{x\m } \tilde{\af }_y (\tilde{\e }_{z})
 \tilde{\e }_{y}  \widetilde{w'}[y,z]) \,  \af _x (1_{x\m } \tilde{\e }_{y z}\m)
 1_{x } \tilde{\af } _x (\tilde{\e }_{y z}) \,  \tilde{\e }_x
\, \widetilde{w'}[x,yz] \tilde{\e }_{xyz}\m = \\ & \af _x (1_{x\m }
\tilde{\af }_y (\tilde{\e }_{z})
 \tilde{\e }_{y})\,  \e _x \, {\e }_x\m  \af _x (1_{x\m }  \widetilde{w'}[y,z]) \,  {\e }_x
\, \widetilde{w'}[x,yz] \tilde{\e }_{xyz}\m = \\  & \af _x (1_{x\m }
\tilde{\af }_y (\tilde{\e }_{z})) \, \af _x (1_{x\m }
 \tilde{\e }_{y})\,  \e _x \,    \widetilde{w'}[x,y] \, \widetilde{w'}[xy,z] \tilde{\e }_{xyz}\m ,
\end{align*} by using the definitions of  $ \tilde{\af } _x $ and
$\e _x,$ and the equality (\ref{cocycle2}). Now, applying
(\ref{composition}) and expanding  $\widetilde {w} [x,y]\m$ according
to  (\ref{tildew}), we come to
\begin{align*} &  \widetilde{w}[x,y] \, \tilde{\af }_{xy} (\tilde{\e }_{z}) \,  \widetilde {w} [x,y]\m
\af _x (1_{x\m }
 \tilde{\e }_{y})\,  \e _x \,    \widetilde{w'}[x,y] \, \widetilde{w'}[xy,z] \tilde{\e
 }_{xyz}\m= \\ & \widetilde{w}[x,y] \, \tilde{\af }_{xy} (\tilde{\e }_{z}) \,
 \tilde{\e }_{xy} \widetilde{w'}[x,y ]\m \tilde{\e }_x\m  \tilde{\af }_x (\tilde{\e}_y)\m
 1_x \, \af _x (1_{x\m }
 \tilde{\e }_{y})\,  \e _x \,    \widetilde{w'}[x,y] \, \widetilde{w'}[xy,z] \tilde{\e
 }_{xyz}\m ,
\end{align*} which, after making the cancellations, gives
\begin{align*}  & 1 _x  \widetilde{w}[x,y] \, \tilde{\af }_{xy} (\tilde{\e }_{z}) \,
 \tilde{\e }_{xy}     \widetilde{w'}[xy,z] \tilde{\e
 }_{xyz}\m  = 1_x \, \widetilde{w}[x,y] \, \widetilde{w}[xy,z],
\end{align*} which is the right hand side of (\ref{extend2}), as desired.

Finally, for a non-transitive  $\af $ we partition the blocks of
$\A $ into  orbits, so that the construction of the $\widetilde{w} $'s
and the verification of (\ref{extend2}), is obviously reduced to
the
transitive case. Our result  follows now from Theorem~\ref{globreduced}.  \fim \\

\end{section}

\begin{section}{Uniqueness}\label{sec:uniqueness}

In this section we extend Definition~\ref{equivalent} permitting
$G$ to act on distinct rings, and prove that any two globalizations of
a twisted partial action of a group $G$ are equivalent,
provided that the rings  are products of blocks.

\begin{defi}\label{iso} Let ${\A }_1$ and  ${\A }_2$ be
rings. We say that a twisted partial action
$${\af }_1 = (\{{\D }^{(1)}_x \}_ {x \in G}, \{{\af}_{1,x} \}_ {x \in G},
 \{w_1[x,y] \}_ {(x,y) \in G\times G})$$  of $G$
on ${\A }_1 $ is isomorphic to the twisted partial action
 $${\af }_2 = (\{{\D}^{(2)}_x \}_ {x \in G}, \{{\af}_{2,x} \}_ {x \in G},
 \{w_2[x,y] \}_ {(x,y) \in G\times G})$$  of $G$ on ${\A }_2$
  if there exist a ring isomorphism
 $\phi : {\A }_1 \to {\A }_2$ such that
 \begin{equation}\label{isoD}{\phi  }( {\D }^{(1)}_x )   = {\D }^{(2)}_x  \, \, \, \forall x \in G,
 \end{equation}
 $${\phi }\, \circ {\af }_{1,x} \circ  {\phi}\m (a) =  {\af }_{2,x} (a) \, \, \,
 \forall x \in G, a \in  {\D}^{(2)}_{x\m}, $$ and
 \begin{equation}\label{isow}{\phi }_{x,y}\, w_1[x,y] \, {\phi }\m _{x,y}   =
 w_2[x,y]  \, \, \, \forall x,y \in G
 \end{equation} as
 multipliers of ${\D }^{(2)}_{x} {\D }^{(2)}_{xy}$ where $ {\phi }_{x,y}$ stands for the restriction of
 $\phi $ to ${\D }^{(1)}_{x} {\D }^{(1)}_{xy}.$ Two globalizations
$${\bt }_1 = ({\B }_1,  \{{\bt}_{1,x} \}_ {x \in G},
 \{u_1[x,y] \}_ {(x,y) \in G\times G})$$ and
$${\bt}_2 = ({\B }_2 ,  \{{\bt}_{2,x} \}_ {x \in G},
 \{u_2[x,y] \}_ {(x,y) \in G\times G})$$ of a partial action
$$\af = (\{\D_x \}_ {x \in G}, \{{\af}_x \}_ {x \in G},
 \{w[x,y] \}_ {(g,h) \in G\times G})$$ of $G$ on $\A $ with embeddings
 ${\f }_1: \A \to {\B }_1$ and ${\f }_2: \A \to {\B }_2$ shall be
 called isomorphic if the actions $\bt _1 $ and $\bt _2$ are isomorphic with the ring isomorphism
 $\phi : \B _ 1 \to \B _2,$ and
 $$ \phi \circ {\f }_1 = {\f }_2 .$$
\end{defi}

We proceed with the next:

\begin{lema}\label{unicitylema1} Let $$\af = (\{\D_x \}_ {x \in G}, \{{\af}_x \}_ {x \in G},
 \{w[x,y] \}_ {(g,h) \in G\times G})$$ be a twisted partial
 action of $G$ on a  unital ring $\A $ and
$${\bt }_1 = ({\B }_1,  \{{\bt}_{1,x} \}_ {x \in G},
 \{u_1[x,y] \}_ {(x,y) \in G\times G}),$$
$${\bt}_2 = ({\B }_2 ,  \{{\bt}_{2,x} \}_ {x \in G},
 \{u_2[x,y] \}_ {(x,y) \in G\times G})$$ be  globalizations of
 $\af $ with embeddings ${\f }_1: \A \to {\B }_1$ and ${\f }_2: \A \to {\B }_2.$
 If
 \begin{equation}\label{sametildew} {\f }\m_1 (u_1[x,y] {\f }_1 (1_{\A })) = {\f }_2\m (u_2[x,y] {\f
 }_2(1_{\A})) \, \, \, \forall x,y \in G,\end{equation}
 then ${\bt }_1 $ and ${\bt
 }_2 $ are isomorphic.
\end{lema}

\p Write $\widetilde{w}[x,y] = {\f }\m_1 (u_1[x,y] {\f }_1 (1_{\A
})),$ $x,y \in G.$ By (ii) of Definition~\ref{globalization}, ${\B
}_1 = \sum_{x \in G} \bt_{1,x}({\f }_1(\A))$ and ${\B }_2 =
\sum_{x \in G} \bt_{2,x}({\f }_2(\A)),$ and we claim that the map
$\phi : {\B }_1 \to {\B }_2$ additively determined by
$\bt_{1,x}({\f }_1(a)) \mt {\bt }_{2,x}({\f }_2(a)),$ $a\in \A ,$
is well defined. For  suppose that $$\sum_{i} {\bt }_{1,x_i} ({\f
}_1(a_i)) =0,$$ and we  want to be sure that $\sum_{i} {\bt
}_{2,x_i} ({\f }_2(a_i)) = 0.$

For all $y \in G$ and $ a \in \A$ we have $\sum_i {\bt }_{1,x_i}
({\f }_1(a_i)) {\bt }_{1,y}({\f }_1(a)) =0$ and applying ${\bt
}\m_y$ we obtain
\begin{align*} & 0 = \sum_i {\bt \m_{1,y}}\circ {\bt
}_{1, x_i} ({\f }_1(a_i)) {\f }_1(a) = \\ & \sum _i u_1[y\m , y]\m
{\bt }_{1,y\m} \circ {\bt }_{1, x_i} ({\f }_1(a_i)) u_1[y\m , y]
{\f }_1 (a) = \\ &  \sum_i  u_1[y\m , y]\m  u_1[y\m , {x_i}] {\bt
}_{1,y\m x_i} ({\f }_1(a_i)) u_1[y\m , x_i]\m u_1[y\m , y] {\f }_1 (a) = \\
 &  \sum_i  u_1[y\m , y]\m  u_1[y\m , x_i] {\bt
}_{1,y\m x_i} ({\f }_1(a_i)) {\f }_1 (1_{y\m x_i})u_1[y\m , x_i]\m
u_1[y\m , y] {\f }_1 (a),
\end{align*} since by (iii) of Definition~\ref{globalization},
${\bt }_{1,y\m x_i} ({\f }_1(a_i)) {\f }_1 (1_{y\m x_i})u_1[y\m ,
x_i]\m u_1[y\m , y] {\f }_1 (a)$  is contained in $${\bt }_{1,y\m
x_i} ({\f }_1(\A ))  {\f }_1 (\A ) = {\f }_1 (\D_{y\m x_i}) = {\f
}_1(\A ) {\f }_1 (1_{y\m x_i}) .$$ Next in view of ${\f }_1(1_{y\m
x_i} ) = {\f }_1 ({\af }_{y\m x_i} (1_{x\m_i y} ) ) = {\bt
}_{1,y\m x_i} ({\f }_1 (1_{x\m_i y}))$ we have that the above sum
equals
\begin{align*}
&  \sum_i  u_1[y\m , y]\m  u_1[y\m , x_i] {\bt }_{1,y\m x_i} ({\f
}_1(a_i 1_{ x\m_i y}))u_1[y\m , x_i]\m u_1[y\m , y] {\f }_1 (a) =
\\
&  \sum_i {\f }_1 ( \widetilde{w}[y\m , y]\m  \widetilde{w}[y\m , x_i]
{\af }_{y\m x_i} (a_i 1_{x\m_i y }) \widetilde{w}[y\m , x_i]\m
\widetilde{w}[y\m , y] a),
\end{align*} which implies
\begin{equation}\label{use} \sum _i  \widetilde{w}[y\m , y]\m  \widetilde{w}[y\m , x_i] {\af
}_{y\m x_i} (a_i 1_{x\m_i y }) \widetilde{w}[y\m , x_i]\m
\widetilde{w}[y\m , y] a = 0,
\end{equation} as ${\f }_1$ is a
monomorphism. Applying ${\f }_2,$ we  obtain by a symmetric
calculation  that
$$\sum_i {\bt }_{2,x_i} ({\f }_2(a_i)) {\bt }_{2,y}({\f }_2(a))
=0$$ for all $y \in G$ and $a\in \A .$ It follows that $\sum_i
{\bt }_{2,x_i} ({\f }_2(a_i)) {\B }_2 =0.$ Now $\B _2$ is a sum of
unital rings and therefore by Remark 2.5 of \cite{DdRS} it is
right  $s$-unital, i.e. for every $b \in {\B }_2$ there exists an
element $e \in {\B }_2 $ such that $b =b e.$ Consequently, $\sum_i
{\bt }_{2,x_i} ({\f }_2(a_i)) =0 ,$ as claimed.

In a similar fashion,   $\bt_{2,x}(\f _2(a)) \mapsto \bt_{1,x}(\f
_1(a)), \;x\in G, a\in \A$ also determines a well-defined map
$\phi': {\B }_2 \to {\B }_1$ and obviously, $\phi' \circ \phi =
\phi \circ \phi' = 1.$ Consequently, $\phi $ is a bijection with
inverse ${\phi }'.$

Observe now that when proving (\ref{use}) we have established that
for all $x, \in G, a,b \in \A ,$
\begin{align*} & {\bt \m_{1,y}}\circ {\bt }_{1, x} ({\f }_1(a)) {\f }_1(b)
=  \f _1 (c)
\end{align*} with $$ c = \widetilde{w}[y\m , y]\m
\widetilde{w}[y\m , x] {\af }_{y\m x} (a 1_{x\m y }) \widetilde{w}[y\m ,
x]\m \widetilde{w}[y\m , y] b,$$   and similarly with ${\bt }_{2,x}$
and ${\f }_2.$ This yields that $\phi $ maps ${\bt _{1,x}} (\f _1
(a)) \cdot {\bt }_{1, y} ({\f }_1(b)) =  {\bt _{1,y}} (  {\bt
\m_{1,y}}\circ {\bt }_{1, x} ({\f }_1(a)) {\f }_1(b)) = {\bt
_{1,y}} ( \f _1 (c) ) $ to  ${\bt _{2,y}} ( \f _1 (c) )  = {\bt
_{2,x}} (\f _2 (a)) \cdot {\bt }_{2, y} ({\f }_2(b)) ,$ so that
$\phi $  preserves multiplication. Since addition is obviously
also preserved,  $\phi $ is a ring isomorphism.

It is immediately  seen that for all $x \in G$ one has $\bt_{2,x}
\circ \phi = \phi \circ \bt_{1,x}$ and $ \phi \circ {\f }_1 = {\f
}_2 .$

It remains to check that $\phi $ respects the twistings, i.e. the
global version of (\ref{isow}). For taking  arbitrary $x, y, z \in
G$ and $a \in \A ,$ we have
\begin{align*} & u_1[x,y] {\bt }_{1,z} ({\f }_1 (a)) =
{\bt }_{1,z} ({\bt }\m_{1,z} (u_1[x,y]) {\f }_1(a) ) = \\ & {\bt
}_{1,z} (u_1[z\m , z]\m {\bt }_{1,z\m} (u_1[x,y])  u_1[z\m ,z]{\f
}_1(a) ) = \\ & {\bt }_{1,z} (u_1[z\m , z]\m u_1[z\m, x] u_1[z\m
x,y] u_1[z\m, xy]\m u_1[z\m ,z]{\f }_1(a) ) =  {\bt }_{1,z} ( {\f
}_1  (d)), \end{align*} with $$ d = \widetilde{w}[z\m , z]\m
\widetilde{w}[z\m, x] \widetilde{w}[z\m x,y] \widetilde{w}[z\m, xy]\m
\widetilde{w}[z\m ,z] a , $$  by the global $2$-cocycle equality and
(\ref{sametildew}). Consequently $\phi $ maps
$$  u_1[x,y] {\bt }_{1,z} ({\f }_1 (a)) = {\bt }_{1,z} (
{\f }_1  (d) ) \mt   {\bt }_{2,z} ( {\f }_1  (d) ) = u_2[x,y] {\bt
}_{2,z} ({\f }_2(a)),$$ by a similar calculation. It follows that
$\phi ( u_1[x,y] b)  = u_2[x,y] {\phi }(b)$ for all $x,y \in G $
and $b \in {\B }_1.$ Analogously, $\phi (b u_1[x,y] )  ={\phi }(b)
u_2[x,y] $ and thus $${\phi }\, u_1[x,y] \, {\phi }\m  =
 u_2[x,y], $$
 for all $x,y \in G , $ which is the global version of (\ref{isow}).
 \fim\\

If $\A$ is a product of blocks, then we are able to show that any two unital globalizations of a partial action $\af $
of $G$ on $\A $ are equivalent in a certain sense, specified below.  For this purpose we will show first that if $\bt $ is globalization for $\af $ with
$\B $ being the ring under the global action, then $\B $ is also a product of blocks, provided that $\B $ has $1_{\B }.$  We first give the details
for the transitive case in Lemma~\ref{transitiveB} below because of simplicity of the notation;
 the general case is obtained similarly which we do in Proposition~\ref{unitalnontransitive}.

Observe next that if $\A $ is a two-sided ideal in a
(non-necessarily unital) ring $\B $ and $\A $ has $1_{\A }$ then
$1_{\A }$ is central in $\B $ and $\A = \B 1_{\A }.$ Moreover, $\A
$ is a direct summand of $\B ,$ for the complement of $\A $ in $\B
$ is the annihilator of $1_{\A } $ in $\B .$ It follows that each
direct summand of $\A $ is a direct summand of $\B .$

\begin{lema}\label{transitiveB} Let $\A $ be a unital ring which is a (finite or infinite) product of
indecomposable rings and let
$${\af } = (\{{\D}_x \}_ {x \in G}, \{ {\af}_{x} \}_ {x \in G},
 \{ w[x,y] \}_ {(x,y) \in G\times G})$$ a transitive  twisted partial action of $G$ on ${\A }$ (in particular,
  each ${\D } _x,$ $x\in G,$ is a unital ring).  Let $H,$  $\Lb $ and ${\Lb }'$  be as in Section~\ref{sec:transitive}, and let
  $${\bt } = ({\B },  \{{\bt}_{x} \}_ {x \in G},
 \{u[x,y] \}_ {(x,y) \in G\times G})$$ be a globalization of $\af ,$ with the embedding $\f : \A \to \B.$ Setting ${\R }_g = \bt _g ( \f (\R _1)), $  $g \in {\Lb
}',$ the following holds:\\

\noindent (i) $\bt _x ({\R }_1 ) = \bt _y ({\R }_1 ) \, \, \, \iff \, \, \, x\m
y \in H.$\\

\noindent (ii)    $\bt _x (\R _g ) = {\R }_{\overline{xg}}$ for each
$g\in {\Lb }'$ and $x \in G.$ \\

\noindent (iii) There is an embedding $$\psi : \B \to {\B }' = \prod _{g\in
{\Lb }'} \R _g $$ which forms commutative triangles with the
projections
 ${\B }' \to {\R }_g $ and $\B \to {\R }_g $
for each $g \in {\Lb }'.$\\

\noindent (iv) $\B $ has $1_{\B }$ if and only if $\psi $ is an isomorphism. In this  case  $\bt $ is transitive.\\

\noindent (v) If $\af $ is unitally globalizable, then any globalization of $\af $ is unital.\\
\end{lema}

\p  Assume for simplicity of notation that $\A \subseteq \B $ so that we can
write ${\R }_g = \bt _g (\R _1),$ with  $g \in {\Lb }'.$ We have
that each block of
 $\A $
 is a direct factor of $\B ,$ and for each $g \in {\Lb }'$ the
 ring ${\R }_g$ is an indecomposable direct factor of $\B ,$ which lies
in $\A $ exactly when $g \in \Lb .$  Evidently $\R _g \neq \R
_{g'} \Longrightarrow \R _g \cap \R _{g'} = 0$ $(g, g' \in {\Lb
}').$ We readily have that $$H = \{x\in G : {\bt}_x ({\R }_1) =
{\R }_1 \}.$$ Next, an invertible multiplier preserves an
idempotent ideal, so that
$$u[x,y] \bt _z (\R _1 ) = \bt _z (\R _1 ) u[x,y] = \bt _z (\R _1 ),$$
for arbitrary $x,y, z \in G.$ In particular, each  $R_g$ $(g \in {\Lb }')$ is preserved by the $u[x,y]$'s.   Now, for $x, y \in G$ we have
$$\bt _x ({\R }_1 ) = \bt _y ({\R }_1 ) \,  \iff \, \R _1 = {\bt }\m _y  \circ  \bt _x (\R _1 ),$$ and the latter equals
\begin{align*} & u[y\m, y]\m  \bt _{y\m} \circ \bt _x (\R _1) \, u[y\m, y] = \\ &
u[y\m, y]\m \, u[y\m,x] \bt _{y\m x} (\R _1)\, u[y\m,x]\m  \, u[y\m, y] =  \bt _{y\m x} (\R _1),
\end{align*} which gives (i).

Next observe that
\begin{align*} \bt _x (\R _1) = u[\overline x , \overline{x}\m x]\m \bt _{\overline x} \circ \bt _{\overline{x}\m x} (\R _1)  u[\overline x , \overline{x}\m x]
= \R _{\overline {x}},
\end{align*} as ${\overline x}\m x \in H.$ Therefore,
\begin{align*} \bt _x (\R _g ) = u[x,g] \bt _{xg} (\R _1 ) u[x,g]\m = \R _{\overline{xg}},
\end{align*} with $x \in G, g \in {\Lb}',$ proving (ii).

  It follows from (ii) that
$$\bt _x (\A ) = \prod _{g\in \Lb } \R _{\overline{xg}} \subseteq
\B ,$$ for any $x\in G.$  Consider the direct product $ {\B }' =
\prod _{g \in {\Lb}'} {\R }_g. $ By the universal property of the
product there is a unique homomorphism $\psi : \B \to {\B }'$
giving commutative triangles with the projections ${\B }' \to {\R
}_g $ and $\B \to {\R }_g $ for each $g \in {\Lb }'.$ Since $\B =
\sum _{x\in G} \bt _x (\A ),$ an arbitrary $b\in \B$ can be
written as $b = \bt _{x_1} (a_1) + \ldots + \bt _{x_s} (a_s)$ for
some $x_1, \ldots x_s \in G$ and $a_1, \ldots, a_s \in \A.$ It is
easy to see that one can choose the elements $a_1, \ldots, a_s \in
\A$ so that if ${\rm pr}_{ g}(  \bt _{x_i} (a_i)) \neq 0$ for some
$i$ and $g\in {\Lb}'$ then ${\rm pr}_{g}(\sum_{j\neq i} \bt _{x_j}
(a_j)) =0,$ where  ${\rm pr}_g$ denotes  the projection $\B \to \R
_g.$ It follows that $\psi $ is injective, proving (iii).

 Now observe
that $\B $ has $1_{\B}$ exactly when there are finitely many $x_1,
\ldots x_s $ in $G$ such that $\B = \sum \bt _{x_i} (\A ).$ Consequently, there exists a partition
\begin{equation}\label{partition}
{\Lb }'= \Lb _1 \cup \Lb _2 \cup \ldots \cup \Lb _s \end{equation} into disjoint subsets
such that \begin{equation}\label{containment}
 \prod _{g\in \Lb _i } \R_g \subseteq \bt _{x_i} (\A ), i=1,\ldots, s, \end{equation} and
\begin{equation}\label{prod}
 \B = (\prod _{g\in \Lb _1} \R _g ) \oplus \ldots   \oplus ( \prod _{g\in \Lb _s} \R _g) \cong
 \prod _{g\in {\Lb}'} \R _g.
\end{equation} This
proves (iv).

Finally, if
\begin{equation}\label{finmany} \{\overline{x_i g} : g\in \Lb ,
i=1, \ldots s \} ={\Lb }',
\end{equation} for some $x_1, \ldots x_s \in G $ then one can obtain a partition (\ref{partition}) with
(\ref{containment}) so that (\ref{prod}) holds. Thus $\B $ has
$1_{\B}$ exactly when (\ref{finmany}) is verified. The latter
condition, however, depends only on $\A , G $ and $\af ,$ and
this proves  (v).    \fim \\

Now we pass to the non-necessarily transitive case. Suppose that
$\A $ is a (finite or infinite) product of blocks:

\begin{equation}\label{decomp3} \A = \prod _{\lb \in
\Lb} {\R}_{\lb},
\end{equation} i.e.  each ${\R}_{\lb} $ is an
indecomposable unital ring, and let
 $${\af } = (\{{\D}_x \}_ {x \in G}, \{ {\af}_{x} \}_ {x \in G},
 \{ w[x,y] \}_ {(x,y) \in G\times G})$$ be a twisted partial action
 of $G$ on $\A $ such that each $\D _x$ $(x\in G)$ is unital. If $\af $ is non-necessarily transitive, then for a given block
  $\R _{\lb}$ its {\it orbit} is defined by
  $  o _{\lb} = \{ {\R }_{{\lb }'} : \exists g \in G , \R _{\lb} \subseteq {\D }_{g\m },
  \af _g (\R _{\lb }) = {\R }_{{\lb }'} \}.$ These are the {\it block-orbits} of $\A $ with respect to
  $\af .$
   Thus we partition the blocks into a disjoint union of  block-orbits,
   and let $\Upsilon \subseteq \Lb $ be such that  $\{\R _{\mu } : \mu \in \Upsilon  \}$ is  a complete
   set of representatives of the
   block-orbits. For any $\mu \in \Upsilon $ put ${\mathcal O} _{\mu } = \prod \R _{{\lb}} $ where
   $\R _{{\lb}}$
   runs over   $o_{\mu}.$  Then evidently
   $$\A = \prod _{\mu \in \Upsilon} {\mathcal O}_{\mu }.$$ The ring ${\mathcal O} _{\mu }$ will be called
   the {\it orbit ideal} corresponding to   ${\mu}.$  Clearly $\af $ restricts to a transitive
   twisted partial action of $G$ on each  orbit ideal.\\

\begin{prop}\label{unitalnontransitive}  Let $\A $ be a unital ring
which is a (finite or infinite) product of indecomposable rings
and let $\af $  be  a unitally  globalizable partial action of $G$
on $\A .$ If ${\bt }$  is an arbitrary  globalization of $\af $ with $G$
acting (globally) on $\B ,$ then\\

 \noindent (i) $\B $ has $1_{\B }$ and $\B $ is a product of blocks;\\

 \noindent (ii) Each block-orbit of $\B $ contains exactly one block orbit of $\A $ and this establishes
 a  one-to-one correspondence between the
 block-orbits of $\B $ and those of $\A ;$\\

 \noindent (iii) The restriction of $\bt $ to an orbit ideal  of $\B $ form a
 globalization for $\af $ restricted to the corresponding
 orbit ideal  of $\A .$
\end{prop}

\p Write $\A $ as in (\ref{decomp3}) and assume the notation
introduced above with respect to the orbits. Assume also for
simplicity  that $\A \subseteq \B .$ Given $\mu \in \Upsilon,$ we have
 a transitive twisted partial
action of $G$ on ${\mathcal O}_{\mu }$ by means of $\af ,$ and let $\R _{\mu}, H(\mu ), {\Lb
}(\mu ), {\Lb }'(\mu )$ correspond to $\R _1, H, \Lb $ and ${\Lb
}',$ respectively, in the notation of Lemma~\ref{transitiveB}.
Similarly as in the transitive case we have an embedding
\begin{equation}\label{embedding}
\psi : \B \to \prod _{\mu \in \Upsilon} \prod_{g\in {\Lb }'(\mu)}
\R _{\mu ,g},\end{equation} where $\R _{\mu, g} = \bt _g (\R _{\mu
}).$ Set ${\mathcal O}'_{\mu } = \sum _{x\in G} \bt _x ({\mathcal
O}_{\mu } ) \subseteq \B .$ Clearly, if $\mu ' \in \Upsilon , \mu
' \neq \mu , $ then $\R _{\mu ' , g} \cap {\mathcal O}'_{\mu } =0$
for any $g \in {\Lb }(\mu '),$  in particular, ${\mathcal O}_{\mu
}$  is the unique orbit ideal of $\A $ contained in $ {\mathcal
O}'_{\mu } .$ It is easily verified that $\bt $ restricts to
${\mathcal O}'_{\mu } $ giving   a globalization for $\af $
restricted to ${\mathcal O}_{\mu } .$ Thus (iii) will follow from
(i) and (ii) provided that the map ${\mathcal O}_{\mu } \mt
{\mathcal O}'_{\mu } $ will give the announced correspondence.

Since $\af $ is unitally globalizable,  we obtain as in Lemma~\ref{transitiveB} that there exist finitely
many elements $x_1, x_2, \ldots , x_s \in G$
such that

\begin{equation}\label{finmany2} \{\overline{x_i g} : g\in \Lb (\mu),
i=1, \ldots s \} ={\Lb }'(\mu),
\end{equation} for all $\mu \in \Upsilon. $ On the other hand, if (\ref{finmany2}) holds, then we see
as in Lemma~\ref{transitiveB} that $\psi $ is an
isomorphism, and in particular, $\B $ has $1_{\B}.$ Again
(\ref{finmany2}) depends only on $\A , G $ and $\af ,$ so that any
globalization of $\af $ is unital. This proves (i). For (ii) it
only remains to note that the ${\mathcal O}'_{\mu } $'s $(\mu \in
\Upsilon)$ are the orbit ideals of $\B $ with
${\mathcal O}'_{\mu } = \prod_{g\in {\Lb }'(\mu)} \R _{\mu ,g}. $ \fim \\

Before stating our main result on uniqueness, we give two more definitions.

\begin{defi}\label{equivalent2} Let
$${\af }_1 = (\{{\D }^{(1)}_x \}_ {x \in G}, \{{\af}_{1,x} \}_ {x \in G},
 \{w_1[x,y] \}_ {(x,y) \in G\times G})$$ be a   twisted partial action of $G$
on a ring ${\A }_1 $ and
 $${\af }_2 = (\{{\D}^{(2)}_x \}_ {x \in G}, \{ {\af}_{2,x} \}_ {x \in G},
 \{ w_2[x,y] \}_ {(x,y) \in G\times G})$$ a  twisted partial action of $G$ on a ring ${\A }_2$ such
 that all
  ${\D }^{(1)} _x, {\D }^{(2)} _x$ $(x\in G)$ are unital rings. We say that ${\af }_1 $ is equivalent
  to ${\af }_2$ if ${\af }_1 $ is isomorphic to a twisted partial
  action of $G$ on ${\A }_2$ which is equivalent to
     ${\af }_2$ in the sense of Definition~\ref{equivalent}.
\end{defi}

The definition   means that  there exist  a ring isomorphism
 $\phi : {\A }_1  \to {\A }_2 $ and   a function $$ G \ni x \mt
{\e}_x \in
 {\U }({\D }^{(2)} _x) \subseteq {\A }_2 $$ such  that (\ref{isoD}) holds,
\begin{equation}\label{alpha2}
{\phi }\, {\af}_{1,x} \, {\phi}\m(a) = {\e}_x \, {\af}_{2,x} (a)
\, {\e}\m _x, \, \, \, \forall x\in G, a \in {\D }^{(2)}_{x\m },
\end{equation} and
 \begin{equation}\label{w2} {\phi } ( w_1[x,y] ) =
 \e _x \,  {\af}_{2,x} (  \e _y 1_{x\m}) \,
 w_2[x,y] \, {\e}\m_{xy}, \, \, \, \forall x,y \in G.
\end{equation} Observe that in this case the $w_1[x,y]$'s are
invertible elements in ${\D }^{(1)}_{x} {\D }^{(1)}_{xy},$ and\\
${\phi }_{x,y}\, w_1[x,y] \, {\phi }\m _{x,y}$ used in
(\ref{isow}) becomes ${\phi } ( w_1[x,y] ).$\\

\begin{defi}\label{equivalentglob}
 Two unital globalizations
$${\bt }_1 = ({\B }_1,  \{{\bt}_{1,x} \}_ {x \in G},
 \{u_1[x,y] \}_ {(x,y) \in G\times G})$$ and
$${\bt}_2 = ({\B }_2 ,  \{{\bt}_{2,x} \}_ {x \in G},
 \{u_2[x,y] \}_ {(x,y) \in G\times G})$$ of a partial action
$$\af = (\{\D_x \}_ {x \in G}, \{{\af}_x \}_ {x \in G},
 \{w[x,y] \}_ {(g,h) \in G\times G})$$ of $G$ on $\A $ with embeddings
 ${\f }_1: \A \to {\B }_1$ and ${\f }_2: \A \to {\B }_2$ shall be
 called equivalent if  ${\bt }_1 $ is equivalent to $\bt _2$ in the sense of
 Definition~\ref{equivalent2} and, moreover,
 $$ \phi \circ {\f }_1 = {\f }_2 ,$$ where $\phi : \B _ 1 \to \B _2$
 is the involved ring isomorphism.
\end{defi}

\begin{teo}\label{uniquenessteo} Let $\A $ be a  ring with $1_{\A }$ which is a (finite or infinite)
product of
indecomposable rings and let
$${\af } = (\{{\D}_x \}_ {x \in G}, \{ {\af}_{x} \}_ {x \in G},
 \{ w[x,y] \}_ {(x,y) \in G\times G})$$ be a unitally  globalizable twisted partial action of $G$ on ${\A }.$
  Then any two
  globalizations of $\af $ are equivalent.
\end{teo}
\p Let $${\bt }_1 = ({\B }_1,  \{{\bt}_{1,x} \}_ {x \in G},
 \{u_1[x,y] \}_ {(x,y) \in G\times G})$$ and
$${\bt}_2 = ({\B }_2 ,  \{{\bt}_{2,x} \}_ {x \in G},
 \{u_2[x,y] \}_ {(x,y) \in G\times G})$$ be arbitrary globalizations
 of $\af $ with embeddings $\f _1 : \A  \to \B _1, \f _2 : \A \to \B _2.$
 By (i) of
 Proposition~\ref{unitalnontransitive},  $\B _1 $ and
 $\B _2$ are unital rings which are products of blocks.  Then clearly each $\B _i$ $(i=1,2)$ is
 the product of the
 block-orbits, and (ii) and (iii) of Proposition~\ref{unitalnontransitive} permit to reduce
 our proof to the case in which  $\af $ is
 transitive. Then by
 Lemma~\ref{transitiveB}, ${\bt }_1$ and $ {\bt }_2$ are also
  transitive. By the global case of
 Proposition~\ref{corestr-equiv}, $\bt _1$ is equivalent to
 $${{\bt }_1}' = ({\B }_1,  \{{\bt}'_{1,x} \}_ {x \in G},
 \{u'_1[x,y] \}_ {(x,y) \in G\times G})$$  and
$${\bt }'_2 = ({\B }_2,  \{{\bt}'_{2,x} \}_ {x \in G},
 \{u'_2[x,y] \}_ {(x,y) \in G\times G})$$
  is equivalent to  ${\bt _2},$ where
$$u'_i[x,y ] =   \prod_{g \in {\Lb }'}  {\bt }\m_{g\m }\circ {\rm pr}_{\f _i (\R _1)}(  u_i[g\m x
\cdot \overline{x\m g}, \, (\overline{x\m g})\m \cdot y \cdot
\overline{y\m x\m g}]),$$  $${\bt }'_{i,x} (a) = {\eta}\m_{i,x}
{\bt }_{i, x} (a) {\eta} _{i,x},$$
$$ \eta _{i,x} =  \prod _{g \in {\Lb }'}  {\bt}\m_{g\m} \circ {\rm pr}_{\f _i (\R _1 )}
(u_i[g\m, x] \, u_i[g\m x \cdot \overline{x\m g}, \,
(\overline{x\m g})\m ]\m) \in {\U}(\B_i),$$ with $x, y \in G,$ $a
\in \B _i ,$  and ${\rm pr}_{\f _i (\R _1 )} $ denoting the
projection $\B _i \to  \f _i (\R _1 ),$ $i=1,2.$ Now, let $\af '$
be the twisted partial action equivalent to $\af $ given by
Proposition~\ref{corestr-equiv}. By the definition of $H,$ for any
$h_1, h_2 \in H$ we have that
$${\rm pr}_{\f _i (\R _1)} u_i[h_1, h_2] = \f _i ( {\rm pr}_{\R
_1} w[h_1, h_2]),$$ and since $g\m x \cdot \overline{x\m g}, \,
(\overline{x\m g})\m \cdot y \cdot \overline{y\m x\m g} \in H,$ it
follows that
\begin{equation}\label{sametildew2} u'_i[x,y] \f _i
(1_{\A } ) = \f _i (\widetilde{w'}[x,y]) \end{equation} and
consequently
$$u'_i[x,y] \f _i (1_x 1_{xy} ) = \f _i ({w'}[x,y]),$$ with
arbitrary $x, y \in G.$

Next, by (ii) of Lemma~\ref{easy} and (\ref{fromeasy}), for $g \in
{\Lb }'$ and $x \in G$ one has
$$ \R _g \subseteq \D _x  \iff  \R _1
\subseteq \D _{g\m x}.$$ Thus if $\R _g \subseteq \D _x ,$ then
$$\R _1 \subseteq \D _{g\m }\D _{g\m x} \ni w[g\m ,x]$$ and $$\R _1
\subseteq \D _{g\m x \cdot \overline{x\m g}} \, \D _{g\m x} \ni
w[g\m x \cdot \overline{x\m g}, \, (\overline{x\m g})\m ] ,$$ so
that $${\rm pr}_{\f _i (\R _1)} u[g\m ,x] = \f _i ( {\rm pr}_{\R
_1} w[g\m ,x]),$$ and $${\rm pr}_{\f _i (\R _1)} u[g\m x \cdot
\overline{x\m g}, \, (\overline{x\m g})\m ] = \f _i ( {\rm pr}_{\R
_1} w[g\m x \cdot \overline{x\m g}, \, (\overline{x\m g})\m ]).$$
This yields that $$ \eta _{i,x} \f _i(1_x 1_{xy }) = \f _i (\e
_x),$$ for any $x \in G.$ Consequently,$$\bt '_{i,x} (\f _i (a)) =
\f _i (\af '_x (a)),$$ for any $a \in \D _{x\m },$ and we conclude
that each $\bt '_i$ $(i=1,2)$ is a globalization of $\af '.$

Next, (\ref{sametildew2}) means that (\ref{sametildew}) of
Lemma~\ref{unicitylema1} is verified, and, consequently, $\bt '_1$
and $\bt '_2$ are isomorphic globalizations of $\af '.$ It follows
that $\bt _1$ and $\bt _2$ are equivalent globalizations of $\af
.$ \fim

\end{section}

\begin{center}
 \large{Acknowledgments}
\end{center}

This work was carried out during the first author's  visits to  the Federal University of Santa Catarina  (Brazil) and the University of  Murcia (Spain). He is grateful to colleagues from both universities for their warm hospitality.

\end{document}